\documentclass[11pt,twoside]{article}
\usepackage{epsfig,amsfonts,color}
\usepackage{amsmath,bm}
\usepackage{mathrsfs}
\allowdisplaybreaks
\usepackage{multirow}
\usepackage{cite}
\usepackage{amssymb, palatino, geometry,url}
\newcommand{\mynote}[1]{{\color{black}{#1}}}
\usepackage[colorlinks=true,linkcolor=blue,citecolor=blue,urlcolor=blue]{hyperref}
\usepackage{lineno}
\usepackage{enumerate}
\usepackage{float}
\usepackage{graphicx}
\usepackage{array}
\usepackage{placeins}
\usepackage{caption}
\usepackage{subcaption}
\usepackage{pdflscape}
\usepackage{algorithm}
\usepackage{algpseudocode}
\usepackage{fancyhdr}
\newcommand{\be}{\begin{equation}}
\newcommand{\ee}{\end{equation}}
\newcommand{\bea}{\begin{eqnarray}}
\newcommand{\eea}{\end{eqnarray}}

\newcommand{\bvec}{\left(\begin{array}{c}}
\newcommand{\evec}{\end{array}\right)}
\newcommand{\bsub}{\begin{subequations}}
\newcommand{\esub}{\end{subequations}}

\newcolumntype{L}[1]{>{\raggedright\let\newline\\\arraybackslash\hspace{0pt}}m{#1}}
\newcolumntype{C}[1]{>{\centering\let\newline\\\arraybackslash\hspace{0pt}}m{#1}}
\newcolumntype{R}[1]{>{\raggedleft\let\newline\\\arraybackslash\hspace{0pt}}m{#1}}

\begin{document}

\title{Stochastic Model Predictive Control for Central HVAC Plants}

\author{Ranjeet Kumar${}^{\P}$, Michael J. Wenzel${}^\ddag$, Mohammad N. ElBsat${}^\ddag$, Michael J. Risbeck${}^\ddag$,\\ Kirk H. Drees${}^\ddag$, Victor M. Zavala${}^{\P}$\\
{\small ${}^{\P}$Department of Chemical and Biological Engineering}\\
{\small \;University of Wisconsin-Madison, 1415 Engineering Dr, Madison, WI 53706, USA}\\
{\small ${}^\ddag$Johnson Controls International}\\
{\small \;507 E. Michigan St., Milwaukee, WI 53202, USA}}
 \date{}
\maketitle

\begin{abstract}
We present a stochastic model predictive control (MPC) framework for central heating, ventilation, and air conditioning (HVAC) plants. The framework uses real data to forecast and quantify uncertainty of disturbances affecting the system over multiple timescales (electrical loads, heating/cooling loads, and energy prices). We conduct detailed closed-loop simulations and systematic benchmarks for the central HVAC plant of a typical university campus. Results demonstrate that deterministic MPC fails to properly capture disturbances and that this translates into economic penalties associated with peak demand charges and constraint violations in thermal storage capacity (overflow and/or depletion). Our results also demonstrate that stochastic MPC provides a more systematic approach to mitigate uncertainties and that this ultimately leads to cost savings of up to 7.5\% and to mitigation of storage constraint violations. Benchmark results also indicate that these savings are close to ideal savings  (9.6\%) obtained under MPC with perfect information. 
\end{abstract}


\section{Introduction} 
Commercial buildings are responsible for over 20\% of the total energy consumption in the U.S. and annual expenditures of over \$200 billion \cite{patel2019}. In this context, heating, ventilation and air-conditioning (HVAC) systems are the largest source of energy use (nearly 50\%) \cite{doe2011}. Central HVAC plants are sophisticated systems that connect multiple energy carriers (water, electricity,  natural gas, cooling water, hot water, steam) and equipment units (pumps, heat exchangers, cooling towers, chillers, and boilers) to meet the cooling and heating loads of single buildings or collections of buildings (e.g., university campuses and urban districts) \cite{rawlings2018}. A central HVAC plant is the equivalent of a utility plant in a manufacturing facility. Besides total energy use, temporal profiles and peak use are also key factors that affect the efficiency and sustainability energy infrastructures. In particular, temporal profiles and peaks might push infrastructures to their design limits (e.g., capacity and ramping) and this might force operators to use inefficient back-up systems. Time-varying market prices and demand charges are used by operators and utility companies to try to mitigate such impacts. These pricing structures create an incentive for HVAC plants to incorporate thermal energy storage (TES) in order to shift loads in time and manipulate peak demands \cite{henze2004,henze2008,ma2012,risbeck2016,risbeck2017}. Effective operation of HVAC plants requires careful real-time management of the multiple components of the central plant; this  is a challenging task because of the tight interconnection of equipment units, the presence of constraints, and the presence of multiple time-varying disturbances (e.g., energy loads and prices). Specifically, disturbances cannot be perfectly anticipated and thus complicate the planning process. All these factors are forcing commercial buildings to incorporate more sophisticated automation systems. 

Model predictive control (MPC) is becoming a established automation technology in HVAC central plants  \cite{avci2013,deng2014,oldewurtel2012,mendoza2012,kwadzogah2013,rawlings2018}. MPC can anticipate and counteract  disturbances and accommodate complex models, constraints, and cost functions \cite{rawlings2018,kumar2019,qin2003survey}. However, existing MPC implementations for HVAC central plants use deterministic representations of the disturbances. In particular, a most likely value (e.g., mean forecast obtained from autoregressive models) is used to compute control actions. Uncertainty associated with forecast errors is thus ignored during the computation of the control action and, instead, errors are counteracted through feedback. This deterministic approach is intuitive and works well in practice but might lead to cost degradation and failure to satisfy constraints \cite{magni2007robustness,limon2010robust}. These issues are often overlooked in the MPC literature because benchmarking procedures often fail to systematically account for the effect of uncertainty (e.g., perfect forecasts are often assumed). 

\mynote{Uncertainty can be explicitly captured in the controller formulations such as stochastic MPC and robust MPC formulations\cite{kerrigan2001robust,bernardini2009scenario,lucia2014handling,de2005stochastic,mesbah2016}. Robust MPC seeks to find optimal control actions to counteract extreme scenarios, whereas stochastic MPC seeks to determine the optimal actions by taking into consideration the probability of all possible occurrences. In a general multi-stage stochastic MPC, one assumes that uncertainty reveals progressively over time (at every stage) and this effect is modeled in the form of a scenario tree. This approach is intuitive as it captures how recourse would be implemented in an ideal setting but the scenario tree grows exponentially with the length of the prediction horizon \cite{lucia2014handling}. The computational intractability of multi-stage MPC is often handled by using a two-stage approximation. Here, it is assumed that all uncertainty reveals after the first time stage and thus recourse is simplified. Some other alternatives for representing recourse include affine decision rules but these approaches tend to decrease flexibility \cite{farina2016,paulson2017}. In this work, we focus on a scenario-based two-stage stochastic MPC for HVAC central plants because we must consider long planning horizons.} 

In the context of energy systems, it has been recently reported that stochastic MPC can systematically mitigate constraint violations and improve economic performance \cite{kumar2019,kumar2018stochastic}. The benefits of stochastic MPC have also been widely reported in the context of building climate (airside) control \cite{oldewurtel2013stochastic,zhang2013scenario,ma2012fast,drgovna2013explicit}  and energy management  \cite{nghiem2017data,garifi2018stochastic,ferrarini2014distributed,ma2012model,ma2014stochastic}. We highlight that these studies have focused on the building (airside); to the best of our knowledge, stochastic MPC formulations for HVAC central plants have not been reported. 

In this work, we present a computational framework for stochastic MPC for HVAC central plants. Our framework addresses HVAC plants for university campuses and seeks to assess the benefits of stochastic MPC over deterministic MPC. The framework uses real disturbance data to conduct forecasting and uncertainty quantification of disturbances. Our benchmarking procedure uses extensive closed-loop simulations under myriad realizations of disturbances in order to properly account for the effect of uncertainty in controller performance. Results indicate that deterministic MPC leads to violations of storage capacity constraints (overflow or drying up) \mynote{of the hot and chilled water tanks} and that stochastic MPC mitigates this issue. We find that storage capacity violations can be partially mitigated in deterministic MPC by adding buffer (back-off) terms but also that stochastic MPC consistently outperforms deterministic MPC in terms of cost. Specifically, we show that stochastic MPC achieves savings in total cost of up to 7.52\%. When these savings are disaggregated, we find consistent reductions in electricity cost (of 6.89\%), in peak demand charges (of 29.8\%), in natural gas cost (of 8.57\%). We also find that stochastic MPC achieves significant reductions in natural gas usage and thus provides an effective approach to manipulate both electricity and natural gas demand profiles. 

\section{Computational Framework}

In this section, we describe the computational framework used in our studies.  We describe the decision-making setting, physical dynamic model, and disturbance forecasting and uncertainty quantification procedures. The framework incorporates deterministic, stochastic, and perfect information MPC formulations. 

\subsection{Decision-Making Setting}

The central HVAC plant for a typical university campus needs to produce chilled water and hot water in order to meet the time-varying loads (demands) from all buildings. The HVAC plant that we consider in this work consists of a chiller subplant that produces chilled water and a heat recovery (HR) chiller subplant that produces both chilled water and hot water, a hot water generator to produce hot water, cooling towers to reduce the temperature of the water purchased from the market, a dump heat exchanger (dump HX) for rejecting heat from the hot water, and storage tanks (one for chilled water and one for hot water). The goal is to determine hourly operating strategies for all equipment units so that the total cost of the external utilities that need to be purchased from the market (electricity, water, and natural gas) is minimized. Water and natural gas are charged on a total demand basis (at time-constant price) while total electricity is charged based on time-varying prices and the monthly peak electrical load is charged based on demand charges.  

The various cost components faced by the central plant are:
\begin{itemize}
\item {\em Electricity transactions (hourly)}: The central plant purchases electricity required by the equipment for their operation. The transactions are charged at the time-varying market price, $\pi^{e}_{t}$. 
\item {\em Water transactions (hourly)}: The central plant needs to purchase water to make up for evaporative losses of water in the cooling towers. Water is purchased from the utility at a fixed price of $\pi^{w}_{t}$= \$0.009/gal. 
\item {\em Natural gas transactions (hourly)}: The central plant needs to purchase natural gas to run the hot water generator to satisfy the campus heating load. Natural gas is available from the utility at a fixed price of $\pi^{ng}_{t}$= \$0.018/kWh. 
\item {\em Peak Electrical Demand Charges (monthly)}: The total electrical load (i.e., the central plant and attached campus load) is charged for the peak demand incurred over a month by the utility company (at a fixed demand charge price of $\pi^{D}$= \$4.5/kW). 
\end{itemize}

\begin{figure}[!htp]
\centering
\includegraphics[width=0.8\textwidth]{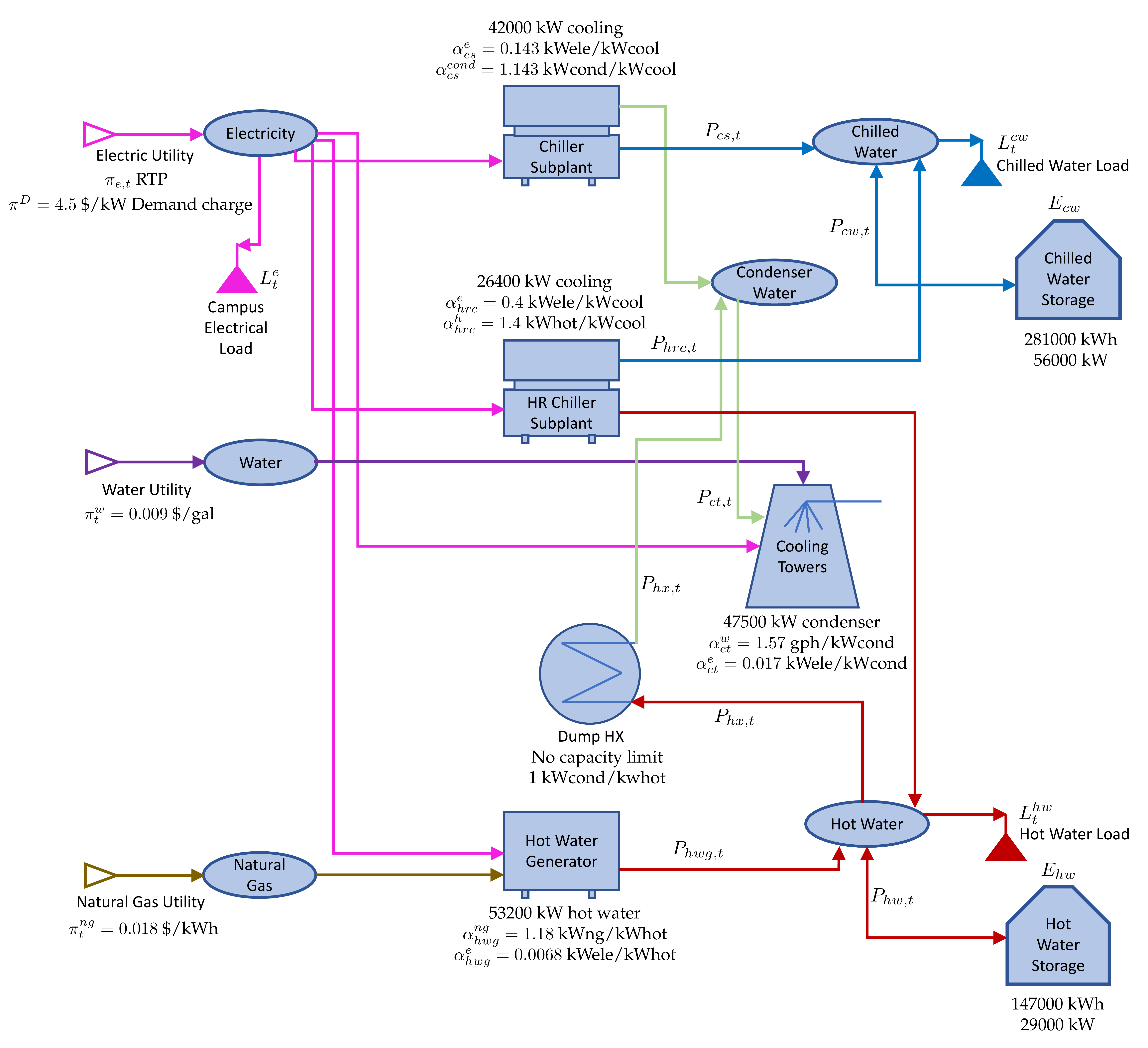}
\caption{Schematic representation of central HVAC plant under study.}
\label{fig:schematic}
\end{figure}

The HVAC central plant considered in this work is illustrated in Figure \ref{fig:schematic}. This shows the energy flows between all the units and interactions with campus loads and utilities. The chiller subplant, HR chiller subplant, hot water generator, and cooling towers consume electricity in their operation. The cooling towers also consume utility water to make up for evaporative losses. The hot water generator is the only unit that consumes natural gas. The electricity, water, and natural gas consumption of these units is tied to their operating loads. The chiller subplant and HR chiller subplant consume $\alpha^{e}_{cs}$ and $\alpha^{e}_{hrc}$ kW of electricity per kW of chilled water produced, respectively; the hot water generator consumes $\alpha^{e}_{hwg}$ kW of electricity and $\alpha^{ng}_{hwg}$ kW of natural gas per kW hot water produced, respectively; and the cooling towers consume $\alpha^{e}_{ct}$ kW of electricity and $\alpha^{w}_{ct}$ utility water per kW of condenser water input, respectively. 

The chilled water load ($L^{cw}_{t}$) of the campus is met by the \mynote{chilled water production from} chiller subplant ($P_{cs,t}$), the HR chiller subplant ($P_{hrc,t}$), and the discharge from chilled water storage ($P_{cw,t}$). The hot water load of the campus is met by \mynote{hot water production from} the HR chiller subplant ($\alpha^{h}_{hrc}P_{hrc,t}$), the hot water generator ($P_{hwg,t}$), and the discharge from the hot water storage ($P_{hw,t}$). The dump heat exchanger (HX) recycles excess hot water ($P_{hx,t}$) in the system by cooling it and producing condenser water which, together with the condenser water produced by the chiller subplant and HR chiller subplants, is cooled further by the cooling towers (total $P_{ct,t}$ condenser water is cooled by the towers). \mynote{The manipulated variables for the system are the operating loads of all units, which include the chilled water production by the chiller and HR chiller subplants, hot water production by the hot water generator, discharge rates from the two storage tanks, the cooling load of the cooling towers and the heat exchange load of the dump HX.}

The HVAC plant operations are driven by uncertain and time-varying disturbances, which are given by the campus loads for electricity ($L^{e}_{t}$), chilled water ($L^{cw}_{t}$), and hot water ($L^{hw}_{t}$), and by the electricity prices ($\pi^{e}_{t}$). The goal of the control (management) system of the plant is to determine operating loads for all units and storage levels to meet campus loads and to minimize the overall plant cost. The appendix provides a detailed description of all the variables and quantities involved.  

\FloatBarrier

\subsection{Forecasting and Uncertainty Quantification}\label{subsec:uq}

The proposed MPC formulations solve discrete-time optimal control problems at every time $t\in \mathbb{Z}_+$. These problems use a forecast (prediction) horizon of length $N$ and denoted by the set $\mathcal{T}:=\{t,t+1,t+2,...,t+N-1\}$. Here, index $t$ corresponds to the time instant $t\cdot h$, where $h$ is the sampling time (assumed to be one hour). We implement zero-order hold for all quantities (these are held constant over the time interval $[(t-1)\cdot h,t\cdot h]$). We implement the MPC schemes under a receding-horizon framework (updated every sampling instant) over a simulation horizon of one year (denoted by the set $\mathcal{Y}:=\{0,...,Y\}$). 

The control actions are the equipment operating loads while the states are the storage levels and other carryover quantities (e.g., peak demand). The decision to determine the operating loads is complicated by the fact that the time-varying electricity prices $\pi_{t}^e$ and the campus electrical load, chilled water load, and hot water load exhibit uncertainty. The loads are denoted as $L_{t} = (L_{t}^e, L_{t}^{cw}, L_{t}^{hw})$. These disturbances are modeled as random variables with realizations indexed by $\xi$. At time $t$, a realization of the disturbances over a forecast horizon $\mathcal{T}$ is denoted as $L_{\mathcal{T}}(\xi),\pi^e_{\mathcal{T}}(\xi)$. The forecast trajectory used in deterministic MPC is a specific realization (usually the one with highest probability) and is denoted as ${\hat{L}_{\mathcal{T}}}$ and ${\hat{\pi}^e_{\mathcal{T}}}$. 
To enable compact representation, we encapsulate all disturbances in the random vector $d_{\mathcal{T}} := (\pi^e_{\mathcal{T}},L_{\mathcal{T}})$. We denote a realization of the disturbance $d_{\mathcal{T}}$ vector as $d_{\mathcal{T}}(\xi)$ and we denote the entire set of disturbance realizations using notation $d_{\mathcal{T}}(\Xi)$.

We assume that the campus electrical load, chilled water load, hot water load, and electricity prices are Gaussian (normal) variables of the form $L_{\mathcal{T}} \sim \mathcal{N}(\hat{L}_{\mathcal{T}},\Sigma_{\mathcal{T}}^L)$, and $\pi^e_{\mathcal{T}} \sim \mathcal{N}(\hat{\pi}^e_{\mathcal{T}},\Sigma_{\mathcal{T}}^e)$. Time correlations in these disturbances are captured by the covariance matrices. \mynote{Gaussian uncertainty descriptions are standard in autoregressive models and we have found them to provide satisfactory results in our studies.} The mean and covariance matrices are updated using a receding-horizon scheme based on historical and real-time data (as it becomes available). We denote the observed (actual) disturbance history of $H$ hours at time $t$ as $d_{\mathcal{H}}$, where $\mathcal{H}:=\{t-H,t-H+1,...,t-1\}$. 

In our implementation, the mean and covariance matrices for $L_{\mathcal{T}}$ and $\pi^e_{\mathcal{T}}$ are obtained using autoregressive (AR) models. Specifically, we use time series models of the form:
\begin{subequations}
    \begin{align}
        & L_t = \sum_{k = 1}^{q} \phi_{k}^L L_{t-k} + c^L + \epsilon_{t}^L\\
        & \pi_t^e = \sum_{k= 1}^{q} \phi_{k}^e \pi^e_{t-k} + c^e + \epsilon_{t}^e,
    \end{align}
\end{subequations}
where $q$ is the order of the model; $\phi_t^L$, $\phi_t^T$, $c^L$ and $c^e$ are coefficients that are learned (estimated) from historical data; and $\epsilon_{t}^L$ and $\epsilon_{t}^e$ are white noise sequences. The mean forecasts $\hat{L}_{\mathcal{T}}$ and $\hat{\pi}^e_{\mathcal{T}}$ (most likely realizations for Gaussian variables) are obtained by using the maximum likelihood estimates of the coefficients. Maximum likelihood procedures also provide estimates of the covariance matrices $\Sigma_{\mathcal{T}}^L $ and $\Sigma_{\mathcal{T}}^e$. Explicit techniques for performing maximum likelihood estimation are provided in \cite{box2015time}. In this work, we use standard procedures provided in the {\tt R} software package.  We highlight that, in this work, we only consider measured disturbances and states. Consequently, no state estimation procedure is needed.

\subsection{Deterministic MPC}\label{subsec:det_mpc}

The deterministic MPC controller uses the mean disturbance forecast $\hat{d}_{\mathcal{T}}$ to find the control policy. This policy minimizes the total cost that is forecast over the prediction horizon $\mathcal{T}$. This is done by solving an optimization problem at time $t \in \mathcal{Y}$. If the prediction horizon $\mathcal{T}$ at time $t$ spans a month, the following optimization problem is solved (with the initial conditions provided by the already known current state values at time $t$, $E_{j,t}, ul_{j,t}, ol_{j,t}, j \in \{cw,hw\}$ and $R_{t}$):

\begin{subequations} \label{Eq:detmpc}
\begin{align}
\min \; & \sum\limits_{k \in \mathcal{T}} \sum\limits_{j=\{e,w,ng\}}\hat{\pi}^j_k r^{j}_{k} + \frac{\pi^D}{\sigma_t} R_{t+1} + \sum\limits_{k \in \mathcal{T}}\sum\limits_{j \in \{cw,hw\}} \rho_{j}(ul_{j,k}+ol_{j,k}). \label{eq:objective}\\
\textrm{s.t. }\; & r^{e}_{k} = \sum_{j \in \{cs,hrc,hwg,ct\}} \alpha^{e}_{j} P_{j,k} + \hat{L}^{e}_{k}, \; k \in \mathcal{T}\label{eq:elec_load}\\
& r^{j}_{k} = \alpha^{j}_{i_j} P_{i_j,k}, \; j \in \{w,ng\}, k \in \mathcal{T}, \; i_{w}=ct, i_{ng}=hwg  \label{eq:water_ng_load} \\
& P_{ct,k} = \alpha^{cond}_{cs} P_{cs,k}+P_{hx,k}, \; k \in \mathcal{T} \label{eq:dumphx} \\
& P_{cs,k}+P_{hrc,k}+P_{cw,k}+S^{un}_{cw,k}-S^{ov}_{cw,k} = \hat{L}^{cw}_{k}, \; k \in \mathcal{T} \label{eq:cw_load} \\
& \alpha^{h}_{hrc} P_{hrc,k}+P_{hwg,k}-P_{hx,k}+P_{hw,k}+S^{un}_{hw,k}-S^{ov}_{hw,k}  = \hat{L}^{hw}_{k}, \; k \in \mathcal{T} \label{eq:hw_load} \\
& E_{j,k+1} = E_{j,k} - P_{j,k}, \; j \in \{cw,hw\}, k \in \mathcal{T} \label{eq:edynamics} \\
& ul_{j,k+1} = ul_{j,k} - S^{m}_{j,k}, \;  m \in \{un, ov\}, j \in \{cw, hw\}, k \in \mathcal{T} \label{eq:uldynamics} \\
& ol_{j,k+1} = ol_{j,k} - S^{m}_{j,k}, \;  m \in \{un, ov\}, j \in \{cw, hw\}, k \in \mathcal{T} \label{eq:oldynamics} \\
& R_{t+1} \geq r^{e}_{k}  \label{eq:peakdemand} \\
& R_{t+1} \geq R_{t}  \label{eq:peakdemand_carry} \\
& \underline{E}_{j,k} \leq E_{j,k} \leq \overline{E}_{j,k}, \; j \in \{cw,hw\}, k \in \mathcal{T} \label{eq:ebound}\\
& \underline{P}_j \leq P_{j,k} \leq \overline{P}_j, \; j \in \{cs,hrc,hwg,ct,hx,cw,hw\}, k \in \mathcal{T} \label{eq:p1bound}\\
& S^{m}_{j,k} \geq 0, \;  m \in \{un, ov\}, j \in \{cw, hw\},  k \in \mathcal{T} \label{eq:Sbound}\\
& ul_{j,k} \geq 0, \; j \in \{cw, hw\}, k \in \mathcal{T}   \label{eq:ulbound}\\
& ol_{j,k} \geq 0, \; j \in \{cw, hw\}, k \in \mathcal{T}  \label{eq:olbound}
\end{align}
\end{subequations}

The constraints \eqref{eq:elec_load}-\eqref{eq:water_ng_load} compute the demands of electricity, water,  and natural gas ($r^{e}_t,r^{w}_t,r^{ng}_t$) that need to be purchased from the utility companies. Constraints \eqref{eq:dumphx}-\eqref{eq:hw_load} impose the energy balance for the condenser water.  Constraints \eqref{eq:cw_load} and \eqref{eq:hw_load} are used to ensure that the chilled and hot water loads are met. Slack variables $S^{m}_{j,k}$, $m \in \{un, ov\}, j \in \{cw, hw\}$ are added in order to maintain feasibility in case of unmet (under-production) or overmet (over-production) of chilled water or hot water. The amount of under-production or over-production of chilled and hot water are carried over using the state variables $ul_{j,k}$ and $ol_{j,k}$, $j \in \{cw, hw\}$. These amounts are computed using constraints \eqref{eq:uldynamics} and \eqref{eq:oldynamics}. The unmet and overmet load state variables are penalized  in the objective function by using the factors $\rho_j$, $j \in \{cw,hw\}$. 

The constraints \eqref{eq:edynamics} describe the dynamics of the state-of charge (SOC) for chilled and hot water TES. Constraint \eqref{eq:peakdemand} ensures that $R_{t+1}$ is the peak residual electricity demand over the horizon $\mathcal{T}$. The parameter $R_t$ is the carryover peak demand (the largest demand seen so far in the month). The peak residual demand is also considered a state of the system and is updated as $R_{t+1} = \max\{R_t,r_{t+1}(\xi)\}$. Constraints \eqref{eq:ebound}-\eqref{eq:olbound} provide bounds on the controls and states. The lower bounds for the operating loads of all units are zero except for the chilled water and hot water storage units (i.e., $\underline{P}_j=0$ for $j\in \{cs,hrc,hwg,ct,hx\}$). For the storage units, the lower bounds of discharge rates correspond to the maximum charging rates, which are negative of the maximum discharging rates (i.e, $\underline{P}_j=-\overline{P}_j, j \in \{cw,hw\}$). \mynote{The demand charge is weighted by a time-varying discounting factor $\sigma_t := \min\{(t_m-t)/N,1\}$ (where $t_m$ denotes the last hour of the month $m$). The discounting factor $\sigma_t$ is used to adjust the demand charge because the prediction horizon $N$ of the MPC formulation is shorter than a month (the period for which the peak electrical demand is charged). Also, it is desired to avoid the occurrence of a peak demand closer to the end of the month. Consequently, the time-varying factor $\sigma_t=\min\{(t_m-t)/N,1\}$ is defined such that the demand charge is penalized  higher when the current time reaches closer to the end of the month than that at the beginning of the month.}

An important consideration when benchmarking the performance of deterministic MPC is that the associated control policy might give rise to constraint violations in the storage levels (when the policy faces actual realizations of loads are observed in the time period $(t,t+1)$). This is because the control policy only satisfies constraints under the mean forecast. To avoid such constraint violations, the bounds on the chilled and hot water TES in \eqref{eq:ebound} are modified to include a buffer capacity \mynote{(a back-off term)}.  A buffer value $\beta\in [0,1]$ forces the storage capacity of chilled and hot water to vary between a fraction $\beta$  and (1-$\beta$) of the maximum capacity (a value of $\beta=0$ gives a standard MPC formulation). \mynote{Such back-off terms are used to tune constraint violations under uncertain disturbances in applications of control and scheduling \cite{grosso2014,koller2018,kumar2018stochastic}.} The buffer is updated after every sampling time in order to prevent the controller from being overly conservative. Specifically, when disturbances push the system outside of the standard constraints, the bounds are adjusted until the optimizer restores feasibility (the following discussion provides more details).

Another important consideration when benchmarking performance is that the demand charge is accounted for based on the peak demand carried over at the end of the month. Consequently, the peak demand needs to be reset at the beginning of each month. To implement this over the simulation horizon $\mathcal{Y}$, we introduce the set $\mathcal{T}_M$ of ending time indices for all months. Index $t_m \in \mathcal{T}_M$ denotes the last hour of the month $m$.  If the prediction horizon $\mathcal{T}$ at the current time $t$ spans two months (i.e., $t < t_m$ and $t+N-1 > t_m$ for some $t_m \in \mathcal{T}_M$), then the optimization problem at time $t$ will consider peak demand variables $R_{1,t+1}$ and $R_{2,t+1}$ for the two months spanned. If at time $t$ in the simulation, $t < t_m$ and  $t+N-1 > t_m$ for some $t_m \in \mathcal{T}_M$, we modify constraints \eqref{eq:peakdemand} and \eqref{eq:peakdemand_carry} as follows:
\begin{subequations}
\begin{align}
& R_{1,t+1} \geq r^{e}_{k}, \; k=t+1, t+2, \dots, t_m \label{eq:mod_peakdemand1} \\
& R_{2,t+1} \geq r^{e}_{k}, \; k=t_m+1, t_m+2, \dots, t+N-1 \label{eq:mod_peakdemand2} \\
& R_{i,t+1} \geq R_{i,t}, \; i = 1, 2.  \label{eq:mod_peakdemand_carry} 
\end{align}
\end{subequations}
where $R_{2,t} = 0$ for $t=t_m$. The objective function \eqref{eq:objective} is also modified to include the  demand charges as:
\begin{align}
\min \; & \sum\limits_{k \in \mathcal{T}} \sum\limits_{j=\{e,w,ng\}}\pi^j_k r^{j}_{k} + \sum_{i=1}^2 \frac{\pi_i^D}{\sigma_t} R_{i,t+1} + \sum\limits_{k \in \mathcal{T}}\sum\limits_{j \in \{cw,hw\}} \rho_{j}(ul_{j,k}+ol_{j,k}). \label{eq:mod_objective}
\end{align}
When $t > t_m+1$ and $t < t_m+2$, the optimization formulation \eqref{eq:objective}-\eqref{eq:Sbound} is solved (with $R_t = R_{2,t_{m+1}}$) because the prediction horizon spans only a single month ($m+1$).

We denote the optimization problem solved in deterministic MPC at time $t$ as $\mathcal{P}^{det}({x_t,\hat{d}_{\mathcal{T}}})$. Here, the input arguments are the current state, control information, and disturbance forecast needed to solve the problem. We define the control action generated from the solution of the problem as $u_t=(P_{j,t},j\in \{cs,hrc,hwg,ct,hx,cw,hw\})$ and the predicted states as $x_{t+1} = (E_{j,t+1},ul_{j,t+1},ol_{j,t+1},R_{t+1}, j \in \{cw,hw\})$ (obtained assuming forecast $\hat{d}_{t+1}$). The control action $u_{t}$ is implemented in the system but the actual disturbance realized $d_{t+1}(\xi)$ will tend to deviate from the forecast $\hat{d}_{t+1}$. As a result, the state at time $t+1$ will differ from that predicted by the MPC controller. To account for the error in the prediction, the horizon is shifted forward to update the disturbance forecast using the new data history. 

After obtaining the solution of $\mathcal{P}^{det}({x_t,\hat{d}_{\mathcal{T}}})$, during the zero-order hold step (inter-optimization period), the current time decisions, $u_t$, for all units are held constant between $(t,t+1)$. During this inter-optimization period, the loads in real time can be varying, and the load balance constraints \eqref{eq:water_ng_load} (production equal to consumption) will not hold. The loads will be balanced by charging or discharging the storage tanks as required in real time. Because the actual discharge does not exactly follow the discharge predicted by the optimization, the states will be at a different level than that predicted by the model at the end of zero-order hold. The amount of the difference is equal to the integrated difference between the constant discharge rate assumed by the optimization model and the actual time-varying discharge rate resulting from load following. We \mynote{propose to} model this integrated difference between the predicted discharge rate and the actual discharge rate by a normal random variable that is added to the predicted energy levels of the TES to simulate an actual energy level in the TES prior to solving the optimization problem on the next time step. The random variable is given by $v_{j,t} = \mathcal{N}(-0.5(L^{j}_{t+1}-L^{j}_{t}), 0.25\sigma^2_{j,err,t+1}+\sigma^2_{j,int})$, $j \in \{cw,hw\}$, where $\sigma^2_{j,err,t}$ is the variance of load prediction error for $t+1$ and $\sigma^2_{j,int}$ is the variance of integrated load for one-hour periods (obtained from historical data). The mean of the random variable $v_{j,t}$ is negative if the load rises in real time during the zero-order hold as the discharge from the TES will increase to make up for the higher load, and similarly, the mean of $v_{j,t}$ is positive if the load falls during the zero-order hold.

We summarize the receding-horizon scheme for deterministic MPC as:
\begin{enumerate}
\item START at $t = 0$ with the given $E_{j,0}$, and, $ul_{j,0}=0, ol_{j,0}=0$ for $j \in \{cw,hw\}$, $R_0=0$. Set $\underline{E}_{j,0} = \beta \overline{E}_j$, $\overline{E}_{j,0} = (1-\beta) \overline{E}_j$, $j \in \{cw,hw\}$. REPEAT for $t\in \mathcal{Y}$:
\item Use disturbance history $d_{\mathcal{H}}$ to obtain forecast $\hat{d}_{\mathcal{T}}$.
\item Solve $\mathcal{P}^{det}(x_t,{\hat{d}_{\mathcal{T}}})$ to obtain decisions $u_{t}=(P_{j,t}, j \in \{cs,hrc,hwg,ct,hx,cw,hw\})$.
\item Implement controls $u_{t}$ over $(t,t+1)$. 
\item Update storage states to actual states as $E_{j,t+1} = E_{j,t}-P_{j,t} + v_{j,t}, j \in \{cw,hw\}$.
\item Obtain the actual realized storage states by subtracting any unmet load from storage states obtained in previous step: $E_{j,t+1} = E_{j,t+1} - ul_{j,t+1}, j \in \{cw,hw\}$. 
\item Update the unmet and overmet load states to $ul_{j,t+1}, ol_{j,t+1}, j \in {cw,hw}$ from Eq. \eqref{eq:uldynamics} and \eqref{eq:oldynamics}. 
\item Modify the bounds on $E_{j,k}$ for $j \in \{cw,hw\}$ in constraints \eqref{eq:ebound} as follows: \\
If $\beta \overline{E}_j \leq E_{j,t+1} \leq (1-\beta) \overline{E}_j$, set $\underline{E}_{j,t+1} = \beta \overline{E}_j$, $\overline{E}_{j,t+1} = (1-\beta) \overline{E}_j$.\\
If $(1-\beta) \overline{E}_j \leq E_{j,t+1} \leq \overline{E}_j$, set  $\underline{E}_{j,t+1} = \beta \overline{E}_j$, $\overline{E}_{j,t+1} = E_{j,t+1}$. \\
If $0 \leq E_{j,t+1} \leq \beta \overline{E}_j$, set  $\underline{E}_{j,t+1} = E_{j,t+1} $, $\overline{E}_{j,t+1} = (1-\beta) \overline{E}_j$. \\
If $E_{j,t+1} \geq \overline{E}_j$,set $E_{j,t+1} =\overline{E}_j$, $\underline{E}_{j,t+1} = \beta \overline{E}_j $, $\overline{E}_{j,t+1} = \overline{E}_j$, and update $ol_{j,k+1}=ol_{j,k+1}+(E_{j,t+1} - \overline{E}_j)$. \\
If $E_{j,t+1} \leq 0$, set $E_{j,t+1} =0$, $\underline{E}_{j,t+1} = 0$, $\overline{E}_{j,t+1} = (1-\beta)\overline{E}_j$, and update $ul_{j,k+1}=ul_{j,k+1}-E_{j,t+1}$.
\item If the current prediction horizon $\mathcal{T}$ spans a single month, update the carry over demand charge $R_{t+1} = \max\{R_{t},r^{e}_{t}\}$, else update as $R_{t+1} = (\max\{R_{t},r^{e}_{t+1}\}, 0)$ with $r^{e}_{t+1}$ being the actual realized residual electrical demand calculated from Eq. \eqref{eq:elec_load} using actual {\em realized} electrical load $L^{e}_{t+1}(\xi)$.
\item Set $t \leftarrow t+1$. If $t = M$, STOP, otherwise RETURN to Step 2.
\end{enumerate} 

In this scheme, the actual realized disturbances are obtained from a set of validation scenarios $\tilde{\Xi}$. These scenarios are generated from actual disturbance data (not from the forecast). The performance of the MPC controller thus depends on the selection of the realized disturbances and is thus random. We run the scheme for the entire set of validation samples to obtain probability distributions for diverse performance metrics such as cost. 

\subsection{Stochastic MPC}\label{subsec:sto_mpc}

In stochastic MPC, uncertainty representations for disturbances are directly captured in the optimization formulation by using multiple realizations (scenarios) $d_{\mathcal{T}}(\xi)$, $\xi\in\bar{\Xi}$, where $\bar{\Xi}$ is a set of sample scenarios. The control decisions $u_{t}=(P_{j,t}, j \in \{cs,hrc,hwg,ct,hx,cw,hw\})$ for the next immediate hour \mynote{(t+1)} are here-and-now (commitment) decisions that need to be made prior to observing uncertainty. The controls $u_k(\xi) = (P_{j,k}(\xi), j \in \{cs,hrc,hwg,ct,hx,cw,hw\})$ at subsequent times $k \in \mathcal{T}\setminus\{t\}$ are modeled as wait-and-see (recourse) variables that can be corrected when the actual disturbance realization is observed. The stochastic MPC problem is thus formulated as a two-stage problem. 

The residual and peak demands are also recourse decisions that are expressed as $r_{\mathcal{T}}(\xi)$ and $\max_{k\in\mathcal{T}}r_k(\xi)$. The SOC of the chilled water and hot water TES at time $t+1$ only depend on the previous storage $E_{j,t}$ and discharge rates $P_{j,t}$, $j \in \{cw,hw\}$. Consequently, $E_{j,t+1}, j \in \{cw,hw\}$ are also here-and-now variables. The rest of the trajectories $E_{j,k}(\xi), j \in \{cw,hw\}$ for $k \in \mathcal{T} \setminus\{t,t+1\}$ are recourse variables because the corresponding $P_{j,k}(\xi), j \in \{cw,hw\}$ are recourse variables for $k \in \mathcal{T}\setminus\{t\}$. 
 
We use $\mathcal{P}^{sto}(x_t,d_{\mathcal{T}}(\bar{\Xi}))$ to denote the optimization problem solved in stochastic MPC at time $t$. The variables and constraints of the formulation are the same as those of the deterministic counterpart but are replicated for the set of realizations $\xi\in \bar{\Xi}$. We use non-anticipativity constraints to enforce the fact that the control actions $u_{t}$ are here-and-now (commitment) decisions that need to be implemented in the system. 

In the proposed stochastic MPC formulation, we solve a stochastic program at time $t \in \mathcal{Y}$ to minimize the expected total forecast cost of the system and satisfy the constraints under all scenarios over the prediction horizon $\mathcal{T}$. If the prediction horizon $\mathcal{T}$ at time $t$ spans a single month, the following problem is solved with the initial conditions provided by the already known current state values at time $t$, $E_{j,t}, ul_{j,t}, ol_{j,t}, j \in \{cw,hw\}$ and $R_{t}$:

\begin{subequations} \label{eq:stochmpc}
\begin{align}
\min \; & \mathbb{E}\left[\sum\limits_{k \in \mathcal{T}} \sum\limits_{j=\{e,w,ng\}}\hat{\pi}^j_k(\xi) r^{j}_{k}(\xi) + \frac{\pi^D}{\sigma_t} R_{t+1}(\xi) + \sum\limits_{k \in \mathcal{T}}\sum\limits_{j \in \{cw,hw\}} \rho_{j}(ul_{j,k}(\xi)+ol_{j,k}(\xi))\right]. \label{eq:objective_sto}\\
\textrm{s.t. }\; & r^{e}_{k}(\xi)= \sum_{j \in \{cs,hrc,hwg,ct\}} \alpha^{e}_{j} P_{j,k}(\xi)+ \hat{L}^{e}_{k}(\xi), \; k \in \mathcal{T}, \xi \in \bar{\Xi} \label{eq:elec_load_sto}\\
& r^{j}_{k}(\xi) = \alpha^{j}_{i_j} P_{i_j,k}(\xi), \; j \in \{w,ng\}, k \in \mathcal{T}, \xi \in \bar{\Xi} \textrm{where $i_{w}=ct, i_{ng}=hwg$} \label{eq:water_ng_load_sto} \\
& P_{ct,k}(\xi) = \alpha^{cond}_{cs} P_{cs,k}(\xi)+P_{hx,k}(\xi), \; k \in \mathcal{T}, \xi \in \bar{\Xi} \label{eq:dumphx_sto} \\
& P_{cs,k}(\xi)+P_{hrc,k}(\xi)+P_{cw,k}(\xi)+S^{un}_{cw,k}(\xi)-S^{ov}_{cw,k}(\xi) = \hat{L}^{cw}_{k}(\xi), \; k \in \mathcal{T}, \xi \in \bar{\Xi} \label{eq:cw_load_sto} \\
& \alpha^{h}_{hrc} P_{hrc,k}(\xi)+P_{hwg,k}(\xi)-P_{hx,k}(\xi)+P_{hw,k}(\xi)+S^{un}_{hw,k}(\xi)-S^{ov}_{hw,k}(\xi) = \hat{L}^{hw}_{k}(\xi), \; k \in \mathcal{T}, \xi \in \bar{\Xi} \label{eq:hw_load_sto} \\
& E_{j,k+1}(\xi) = E_{j,k} - P_{j,k}(\xi), \; j \in \{cw,hw\}, k \in \mathcal{T}, \xi \in \bar{\Xi} \label{eq:edynamics_sto} \\
& ul_{j,k+1}(\xi) = ul_{j,k} - S^{m}_{j,k}(\xi), \;  m \in \{un, ov\}, j \in \{cw, hw\}, k \in \mathcal{T}, \xi \in \bar{\Xi} \label{eq:uldynamics_sto} \\
& ol_{j,k+1}(\xi) = ol_{j,k} - S^{m}_{j,k}(\xi), \;  m \in \{un, ov\}, j \in \{cw, hw\}, k \in \mathcal{T}, \xi \in \bar{\Xi} \label{eq:oldynamics_sto} \\
& R_{t+1}(\xi) \geq r^{e}_{k}(\xi), \; \xi \in \bar{\Xi}  \label{eq:peakdemand_sto} \\
& R_{t+1}(\xi) \geq R_{t}, \; \xi \in \bar{\Xi} \label{eq:peakdemand_carry_sto} \\
& \underline{E}_{j,k} \leq E_{j,k}(\xi) \leq \overline{E}_{j,k}, \; j \in \{cw,hw\}, k \in \mathcal{T}, \xi \in \bar{\Xi} \label{eq:ebound_sto}\\
& \underline{P}_j \leq P_{j,k}(\xi) \leq \overline{P}_j, \; j \in \{cs,hrc,hwg,ct,hx,cw,hw\}, k \in \mathcal{T}, \xi \in \bar{\Xi} \label{eq:p1bound_sto}\\
& S^{m}_{j,k}(\xi) \geq 0, \;  m \in \{un, ov\}, j \in \{cw, hw\},  k \in \mathcal{T}, \xi \in \bar{\Xi} \label{eq:Sbound_sto}\\
& ul_{j,k}(\xi) \geq 0, \; j \in \{cw, hw\}, k \in \mathcal{T}, \xi \in \bar{\Xi}   \label{eq:ulbound_sto}\\
& ol_{j,k}(\xi) \geq 0, \; j \in \{cw, hw\}, k \in \mathcal{T}, \xi \in \bar{\Xi}  \label{eq:olbound_sto}\\
& P_{j,t}(\xi) = P_{j,t}(\xi'), \; j \in \{cs,hrc,hwg,ct,hx,cw,hw\}, \xi \neq \xi', \xi \in \bar{\Xi} \label{eq:nonant_sto}
\end{align}
\end{subequations}

The expected value $\mathbb{E}[\cdot]$ is defined over the set of scenarios $\bar{\Xi}$. We implement the same approach as described for deterministic MPC for resetting the peak electrical demand charge at the beginning of each month, and for updating the bounds on the chilled water and hot water TES in every step of the MPC simulation. For consistency, here we present a stochastic MPC formulation that also accounts for a storage buffer $\beta$. In Section \ref{sec:case} we will see that stochastic MPC with no buffer ($\beta=0$) can outperform deterministic MPC and significantly reduce constraint violations. This is because the controller naturally backs-off from the constraints when multiple disturbance realizations are accounted for. 

We summarize the stochastic MPC scheme as:
\begin{enumerate}
\item START at $t = 0$ with the given $E_{j,0}$, and, $ul_{j,0}=0, ol_{j,0}=0$ for $j \in \{cw,hw\}$, $R_0=0$. Set $\underline{E}_{j,0} = \beta \overline{E}_j$, $\overline{E}_{j,0} = (1-\beta) \overline{E}_j$, $j \in \{cw,hw\}$. REPEAT for $t\in \mathcal{Y}$:
\item Use disturbance history $d_{\mathcal{H}}$ to obtain forecast $d_{\mathcal{T}}(\bar{\Xi})$.
\item Solve $\mathcal{P}^{sto}(x_t,{d_{\mathcal{T}}}(\bar{\Xi}))$ to obtain decisions $u_{t}=(P_{j,t}, j \in \{cs,hrc,hwg,ct,hx,cw,hw\})$.
\item Implement the decisions, $u_{t}$ over $(t,t+1)$. 
\item Update the storage states to the actual states as $E_{j,t+1} = E_{j,t}-P_{j,t} + v_{j,t}, j \in \{cw,hw\}$, where $v_{j,t}$ is the random variable as defined in Section \ref{subsec:det_mpc}.
\item Obtain actual realized storage states by subtracting any unmet load from storage states obtained in previous step: $E_{j,t+1} = E_{j,t+1} - ul_{j,t+1}, j \in \{cw,hw\}$. 
\item Update the unmet and overmet load states to $ul_{j,t+1}, ol_{j,t+1}, j \in \{cw,hw\}$ from Eq. \eqref{eq:uldynamics_sto} and \eqref{eq:oldynamics_sto}. 
\item Modify the bounds on $E_{j,k}$ for $j \in \{cw,hw\}$ in constraints \eqref{eq:ebound_sto} as follows: \\
If $\beta \overline{E}_j \leq E_{j,t+1} \leq (1-\beta) \overline{E}_j$, set $\underline{E}_{j,t+1} = \beta \overline{E}_j$, $\overline{E}_{j,t+1} = (1-\beta) \overline{E}_j$.\\
If $(1-\beta) \overline{E}_j \leq E_{j,t+1} \leq \overline{E}_j$, set  $\underline{E}_{j,t+1} = \beta \overline{E}_j$, $\overline{E}_{j,t+1} = E_{j,t+1}$. \\
If $0 \leq E_{j,t+1} \leq \beta \overline{E}_j$, set  $\underline{E}_{j,t+1} = E_{j,t+1} $, $\overline{E}_{j,t+1} = (1-\beta) \overline{E}_j$. \\
If $E_{j,t+1} \geq \overline{E}_j$, set $E_{j,t+1} =\overline{E}_j$, $\underline{E}_{j,t+1} = \beta \overline{E}_j $, $\overline{E}_{j,t+1} = \overline{E}_j$, and update $ol_{j,k+1}=ol_{j,k+1}+(E_{j,t+1} - \overline{E}_j)$. \\
If $E_{j,t+1} \leq 0$, set $E_{j,t+1} =0$, $\underline{E}_{j,t+1} = 0$, $\overline{E}_{j,t+1} = (1-\beta)\overline{E}_j$, and update $ul_{j,k+1}=ul_{j,k+1}-E_{j,t+1}$.
\item If the current prediction horizon $\mathcal{T}$ spans a single month, update the carry over demand charge $R_{t+1} = \max\{R_{t},r^{e}_{t}\}$, else update as $R_{t+1} = (\max\{R_{t},r^{e}_{t+1}\}, 0)$ with $r^{e}_{t+1}$ being the actual realized residual electrical demand calculated from Eq. \eqref{eq:elec_load_sto} using actual {\em realized} electrical load $L^{e}_{t+1}(\xi)$.
\item Set $t \leftarrow t+1$. If $t = M$ STOP, otherwise RETURN to Step 2.
\end{enumerate} 

In this scheme, the actual realized disturbances are obtained from the set of validation samples $\tilde{\Xi}$. Importantly, these validation samples differ from the realizations used in the MPC controller formulation $\bar{\Xi}$. By running the stochastic MPC scheme for all validation samples, we can compute probability distributions for performance metrics that are compared with those from deterministic MPC. This systematic procedure ensures fair comparisons between different MPC implementations.

\subsection{Perfect Information MPC}\label{subsec:perfinfo}
We also consider a perfect information MPC implementation to evaluate the ideal performance of MPC. Under perfect information MPC, we compute commitment policies $u_t(\xi)$ at every time $t$ for every realization $\xi\in \tilde{\Xi}$ of the loads and prices. These policies can be computed by removing the nonanticipativity constraints \eqref{eq:nonant_sto} from the stochastic MPC formulation \eqref{eq:stochmpc} and replacing the disturbance forecast model with the true disturbance signals corresponding to each realization $d(\xi)$, and implementing the same scheme as the stochastic MPC run over a year period $\mathcal{Y}$ \cite{kumar2018stochastic}. 

\subsection{Handling Constraint Violations}

The deterministic MPC controller can violate the constraints during the transition $(t,t+1)$ if the forecasts are poor or if the scenarios used. The stochastic MPC formulation can also incur violations if the scenarios used do not capture the actual realizations. To capture this issue in our closed-loop simulations, an auxiliary MPC controller is used to correct the control actions and restore feasibility. At time $t$, this controller solves the feasibility restoration problem: 

\begin{subequations}\label{eq:restoration}
    \begin{align}
    \min_{\substack{\Delta P_j, \\ j \in \{cs,hrc,hwg,ct,hx,cw,hw\}}} \; &  \sum_{j}|\Delta P_j| \\
	\textrm{s.t.} \; & E_{j,t+1} = E_{j,t} - P_{j,t} + \Delta P_j, \; j \in \{cw,hw\}\\
	& \sum_{j \in \{cs,hrc,cw\}} (P_{j,t}+\Delta P_{j,t}) = L^{cw}_{t} \\
    & \sum_{j \in \{hwg,hw\}}(P_{j,t}  + \Delta P_{j}) + \alpha^{h}_{hrc}\,(P_{hrc,t}+\Delta P_{hrc})-P_{hx,t}-\Delta P_{hx} = L^{hw}_{t} \\
	& 0 \leq E_{j,t+1} \leq \overline{E}_j, \; j \in \{cw,hw\}\\
	& \underline{P}_j \leq P_{j,t} + \Delta P_j \leq \overline{P}_j, \; j \in \{cs,hrc,hwg,ct,hx,cw,hw\}.
    \end{align}
\end{subequations}
This formulation uses the actual realizations for the chilled water and hot water loads over the time interval $(t,t+1)$. The feasibility restoration problem seeks to find a net rate correction $\Delta P_j, j \in \{cs,hrc,hwg,ct,hx,cw,hw\}$ that satisfies the storage constraints. If the auxiliary controller fails to remain feasible even after solving the restoration problem, we assume that no action is taken. In other words, we set $P_{j,t}=0$ for $j \in \{cs,hrc,hwg,ct,hx,cw,hw\}$ and correct the states $E_{j,t+1}$ for $j \in \{cw,hw\}$ and $R_{t+1}$ accordingly. This leads to loss of performance because there is a failure to meet the campus loads during time $(t,t+1)$.

\subsection{Benchmarking Procedure} \label{subsec:benchmark}
We extend the benchmarking procedure given in \cite{kumar2019} for battery management systems. To distinguish the policies obtained from the three MPC schemes, we denote the policies obtained from the deterministic MPC as $u_t^{det}$ and $r_t^{det}(\xi)=(r^{j,det}_t, j \in \{e,w,ng\})$, those from stochastic MPC scheme as $u_t^{sto}$ and $r_t^{sto}(\xi)=(r^{j,sto}_t, j \in \{e,w,ng\})$, and those from perfect information MPC as  $u_t^{perf}(\xi)$ and $r_t^{perf}(\xi)=(r^{j,perf}_t, j \in \{e,w,ng\})$ (corresponding to each realization $\xi$).

Each realization in the validation set $\tilde{\Xi}$ generates an annual cost for the MPC controllers. The annual cost under stochastic MPC for a given realization $\xi\in\tilde{\Xi}$ and under a given closed-loop policy ${u_{\mathcal{M}},r_{\mathcal{M}}}$ is given by:
\begin{align}\label{eq:costt}
\Phi^{sto}(\xi) := \sum\limits_{t \in \mathcal{Y}} \sum\limits_{j=\{e,w,ng\}}\pi^j_k(\xi) r^{j,sto}_{k}(\xi) +\sum_{m=1}^{12} {\pi^D} \max_{t \in \{1,\dots, t_m\}} r_t^{e,sto}(\xi).
\end{align}

The annual cost for deterministic MPC is denoted as $\Phi^{det}(\xi)$ and is defined as in \eqref{eq:costt} (but with its corresponding closed-loop policy $u_t^{det}$ and $r_t^{det}(\xi)=(r^{j,det}_t, j \in \{e,w,ng\})$). The annual cost for perfect information MPC is:
\begin{align}
\Phi^{perf}(\xi) := \sum\limits_{t \in \mathcal{Y}} \sum\limits_{j=\{e,w,ng\}}\pi^j_k(\xi) r^{j,perf}_{k}(\xi) +\sum_{m=1}^{12} {\pi^D} \max_{t \in \{1,\dots, t_m\}} r_t^{e,perf}(\xi).
\end{align}
The costs for the different validation realizations are used to create empirical probability distributions and cumulative probability distributions for diverse quantities of interest and to compute statistics such as expected costs $\mathbb{E}[\Phi^{det}(\tilde{\Xi})]$, $\mathbb{E}[\Phi^{sto}(\tilde{\Xi})]$, and $\mathbb{E}[\Phi^{perf}(\tilde{\Xi})]$.

Of particular interest in our benchmark studies is a metric that we call the {\em expected cost of the HVAC central plant}. To compute this value, we evaluate the total cost under the assumption that there is no HVAC central plant serving campus. For a particular validation realization $\xi\in \tilde{\Xi}$, this cost is denoted as $\Phi^{nocp}(\xi)$. The ideal expected cost of the central plant is defined as $\textrm{CCP}^{perf}(\xi) := \Phi^{perf}(\xi)-\Phi^{nocp}(\xi)$. The cost of the central plant under stochastic MPC is $\textrm{CCP}^{sto}(\xi) := \Phi^{sto}(\xi)-\Phi^{nocp}(\xi)$ and under deterministic MPC is $\textrm{CCP}^{det}(\xi) := \Phi^{det}(\xi)-\Phi^{nocp}(\xi)$. As in the case of cost, the realizations are used to obtain probability distributions and to compute statistics such as $\mathbb{E}[\textrm{CCP}^{perf}(\tilde{\Xi})]$, $\mathbb{E}[\textrm{CCP}^{sto}(\tilde{\Xi})]$, and $\mathbb{E}[\textrm{CCP}^{det}(\tilde{\Xi})]$. The cost of the central plant is a metric that reflects losses/gains in asset value due to the use of better control policies (it isolates the effect of the control from that of the equipment). We also consider the value of the stochastic MPC, which is defined as $\textrm{VSMPC}(\xi) := \textrm{CCP}^{det}(\xi)-\textrm{CCP}^{sto}(\xi)$, and the expected value of stochastic MPC as $\mathbb{E}[\textrm{VSMPC}(\tilde{\Xi})]$. 

\section{Benchmark Results}\label{sec:case}

We now present closed-loop simulation results for deterministic, stochastic, and perfect information MPC for an entire year of operation. The controllers use a prediction horizon of 168 hours (seven days), which is chosen based on the observation that the data for the loads and electricity prices exhibit weekly periodicity (e.g., high load in weekdays and low load in weekends). In other words, a horizon of seven days captures periodic effects \cite{kumar2018stochastic}. The stochastic MPC problem contains 100 forecast scenarios and a total of 168,450 variables and 143,750 constraints (the realizations are obtained using Monte Carlo sampling). A total of 200 closed-loop year-long runs were performed for each MPC implementation (using the same validation scenarios). The number of forecast and validation scenarios are chosen considering a trade-off between performance and computational cost. \mynote{Specifically, using few scenarios leads to fast ccomputations but the full uncertainty space is not well represented and this limits the benefits of stochastic MPC. Having a larger number of scenarios increases the computational time but the controller performance improves. Hence, we use a heuristic-based approach to chose a suitable number of forecast and validation scenarios to benchmark stochastic MPC.} The simulations were run on a 32-core machine with Ubuntu 14.04, Intel Xeon 2.30 GHz processors and 188 GB RAM. The schemes are implemented in {\tt Julia} and leverage the algebraic modeling capabilities of {\tt JuMP} \cite{dunning2017jump}. The optimization problems are solved in extensive form using {\tt Gurobi 8.1}. 

We use time series forecasting procedures provided in the {\tt R} software package. Specifically, we use the {\tt ar} function in  {\tt R} to estimate the coefficients of the AR model of order $q$ with the settings as {\tt aic=false}, {\tt order=$q$}, and {\tt method="ols"} (ordinary least-squares). The forecasts are obtained using the {\tt forecast} function with prediction horizon of $N$, and covariance matrices for the forecasts are obtained using {\tt var.pred} function. With this information, we generate a set of disturbance realizations $\bar{\Xi}$ by sampling from the corresponding estimated probability density functions.   Forecasts and scenarios for the electrical load, chilled water load, hot water load, and the electricity price are obtained using a $q=$168 order AR model and predicts the loads over an horizon of $N=$168 hours. The parameters of the AR models are estimated by using a data history of 184 days (i.e., 6 months) at every time step in the closed loop of MPC. 

Each MPC problem instance takes\mynote{, on average,} about one second to solve for deterministic MPC and about 6-7 seconds to solve for stochastic MPC. Despite these fast solutions, we note that year-long closed-loop simulations required approximately {\em 2 hours} for deterministic MPC and {\em 10 hours} for stochastic MPC. These computational times comprise forecasting, optimization solution, and feasibility checks. These computational times are serial but can be partially parallelized (this is left as a topic of future work since the computational workflows involved are complex). 

\subsection{Forecasting of Loads and Electricity Prices}\label{subsec:forecasting}

\begin{figure}[!htp]
\centering
\includegraphics[width=0.7\textwidth]{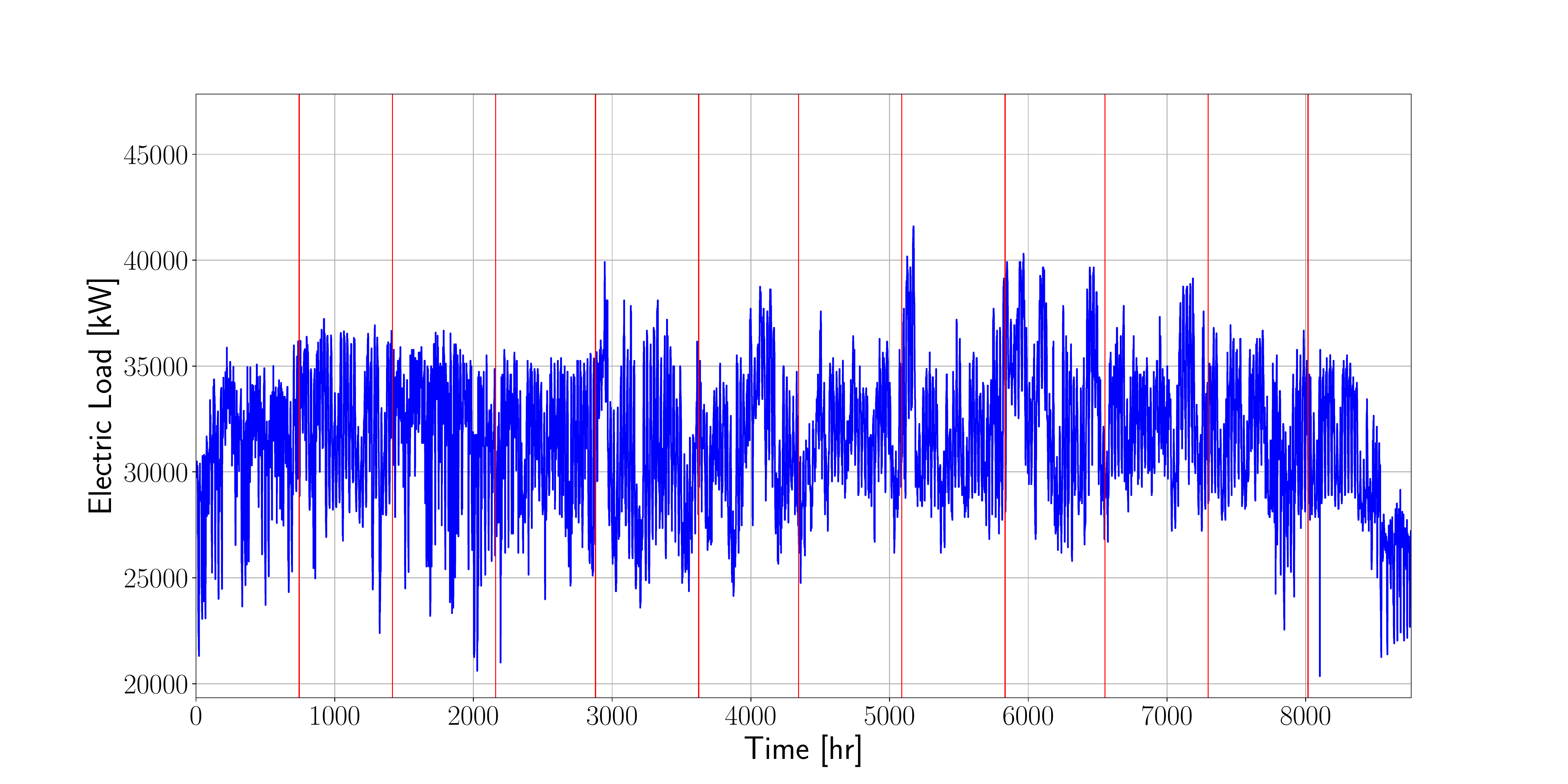}
\caption{Historical electrical load of the campus. Red vertical lines denote end of each month.}
\label{fig:eload}
\end{figure}

\begin{figure}[!htp]
\centering
\includegraphics[width=0.7\textwidth]{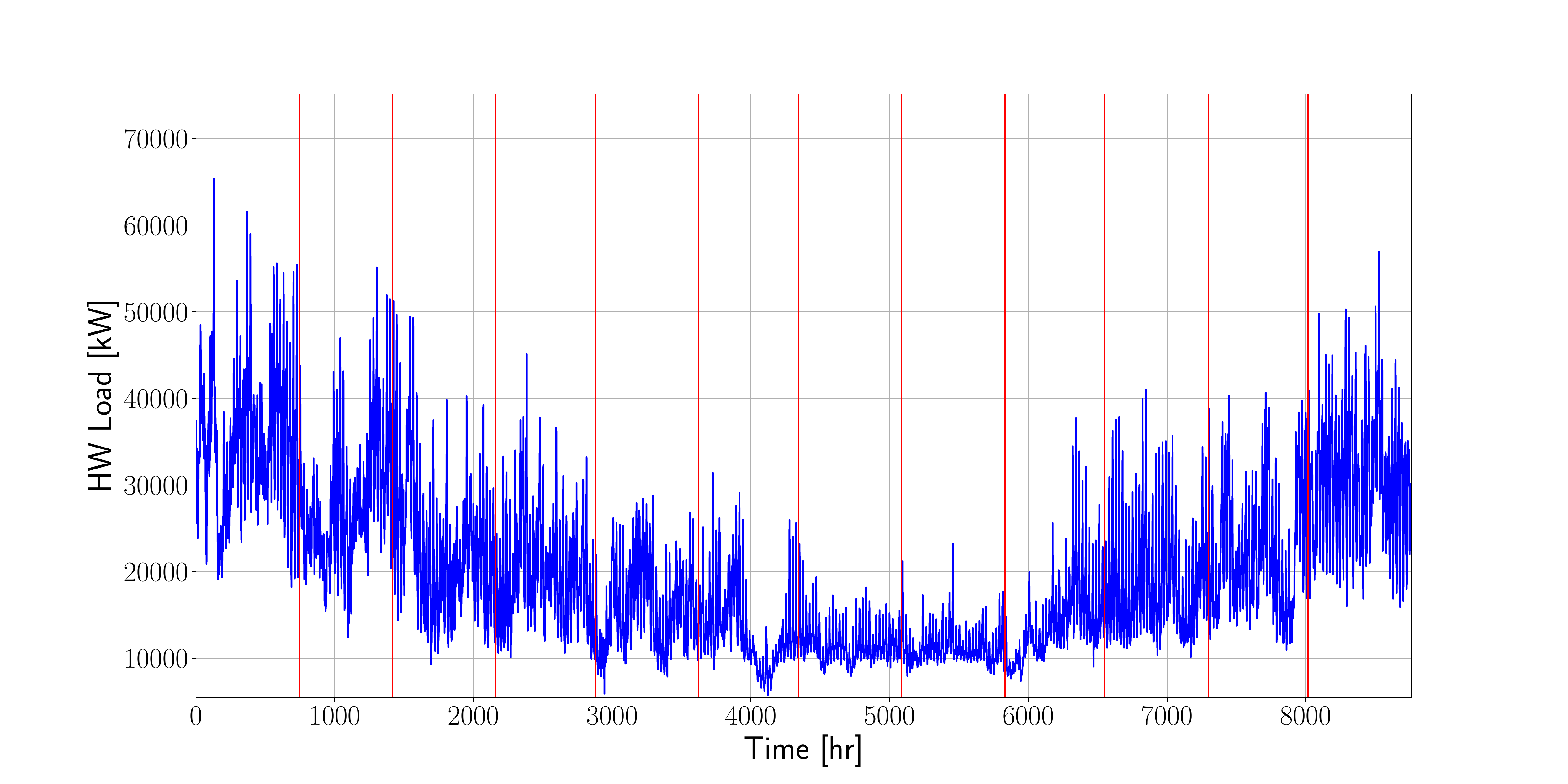}
\caption{Historical hot water load of the campus. Red vertical lines denote end of each month.}
\label{fig:hload}
\end{figure}

\begin{figure}[!htp]
\centering
\includegraphics[width=0.7\textwidth]{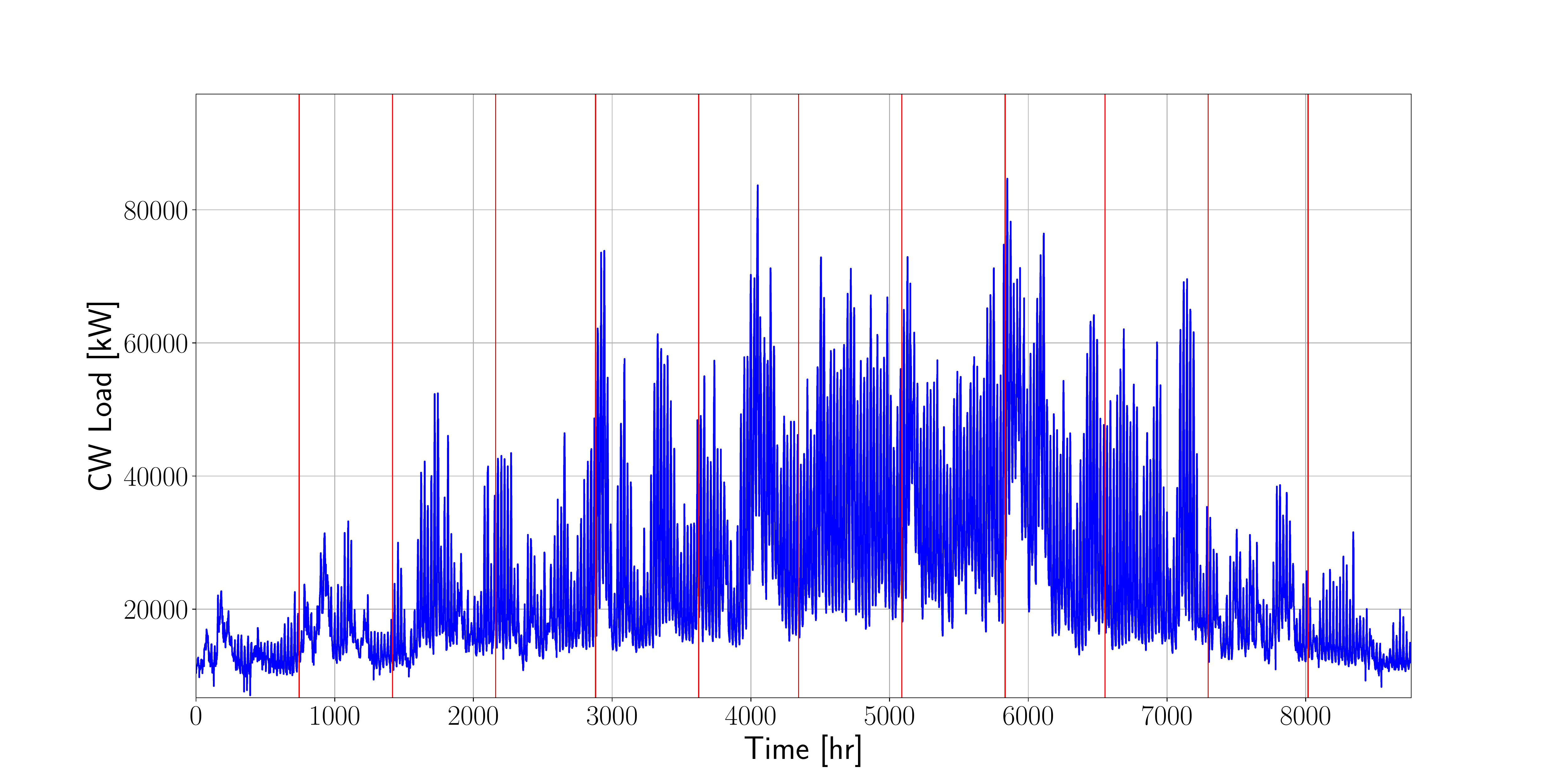}
\caption{Historical chilled water load of the campus. Red vertical lines denote end of each month.}
\label{fig:cload}
\end{figure}

\begin{figure}[!htb]
\centering
\includegraphics[width=0.7\textwidth]{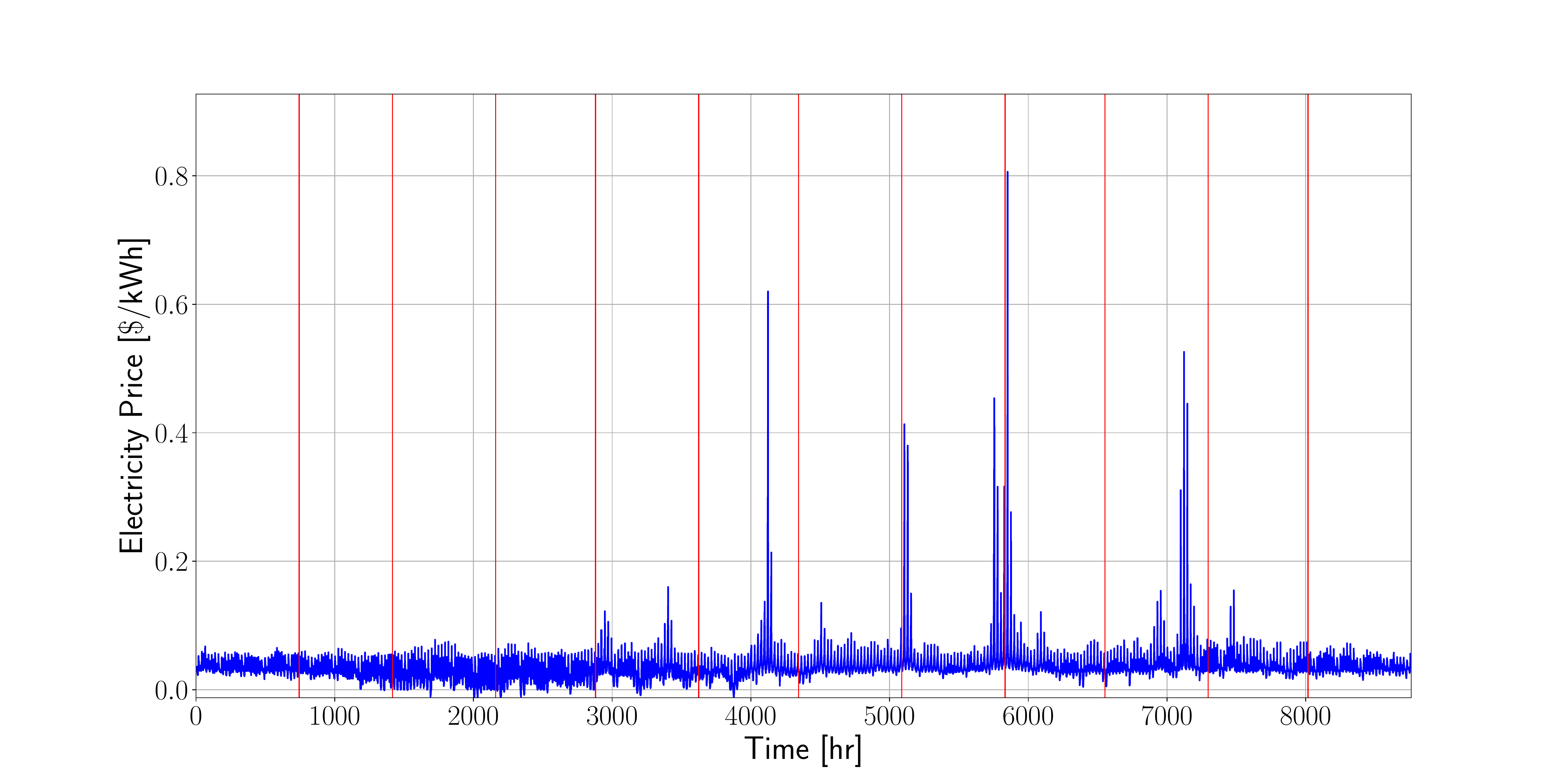}
\caption{Historical electricity price data. Red vertical lines denote end of each month.}
\label{fig:eprice}
\end{figure}

Figures \ref{fig:eload}-\ref{fig:cload} show historical data for the campus electrical load, hot water load, and chilled water load for the entire year.  Figure \ref{fig:eprice} shows historical electricity prices for the same period. The vertical red lines represent monthly periods.

\FloatBarrier
\FloatBarrier
\FloatBarrier

\FloatBarrier

\begin{figure}[!htp]
\centering
\includegraphics[width=0.6\textwidth]{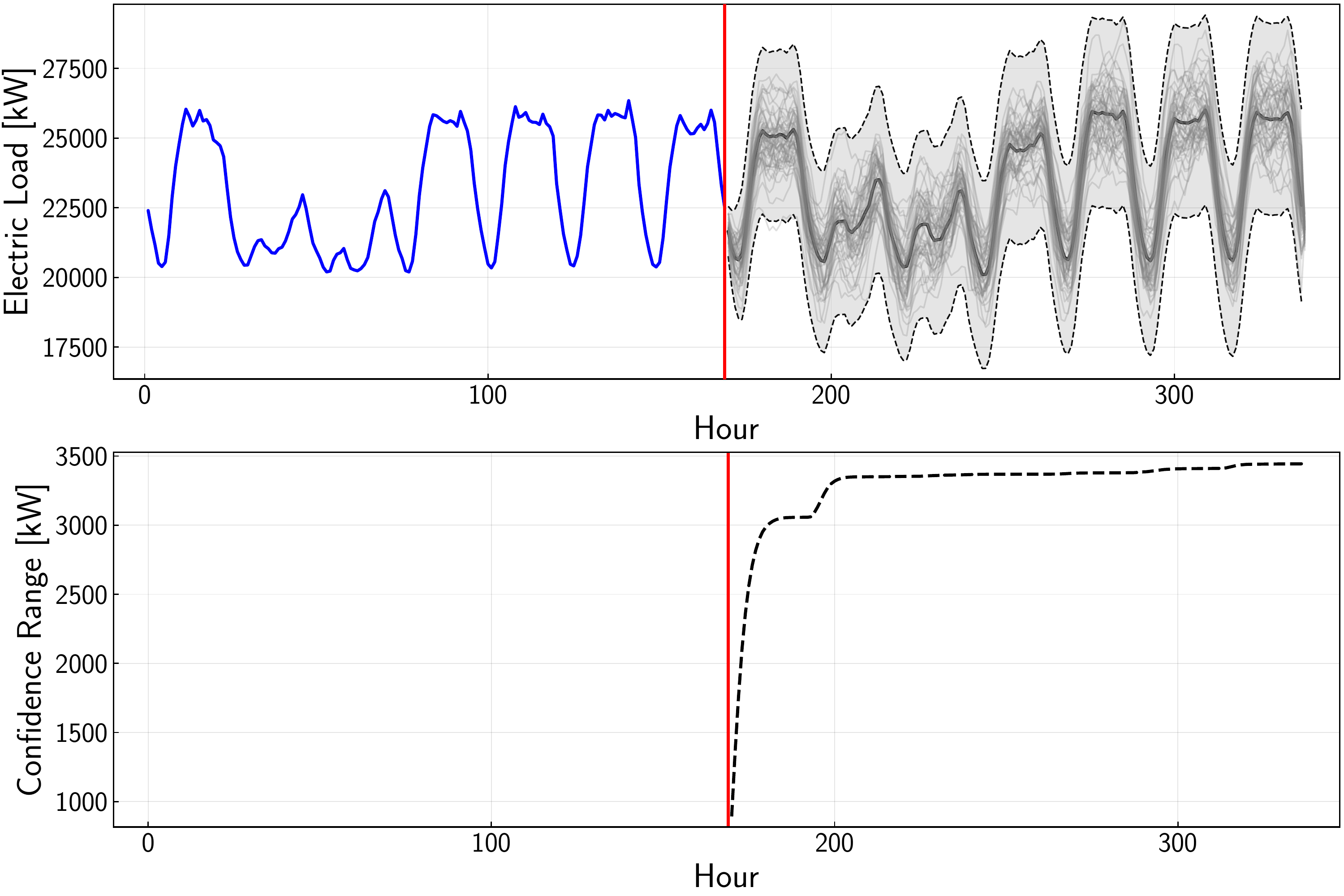}
\caption{A single instance of electrical load forecast for 1-week with AR model. Vertical red line denotes the current time. In the top panel, dark bold curve represents mean forecast, the grey band denotes 99\% confidence interval and the light black curves within the band represent a few sample scenarios. The bottom panel shows the trajectory of the 99\% confidence range with prediction time.}
\label{fig:load_fore}
\end{figure}

\begin{figure}[!htp]
\centering
\includegraphics[width=0.6\textwidth]{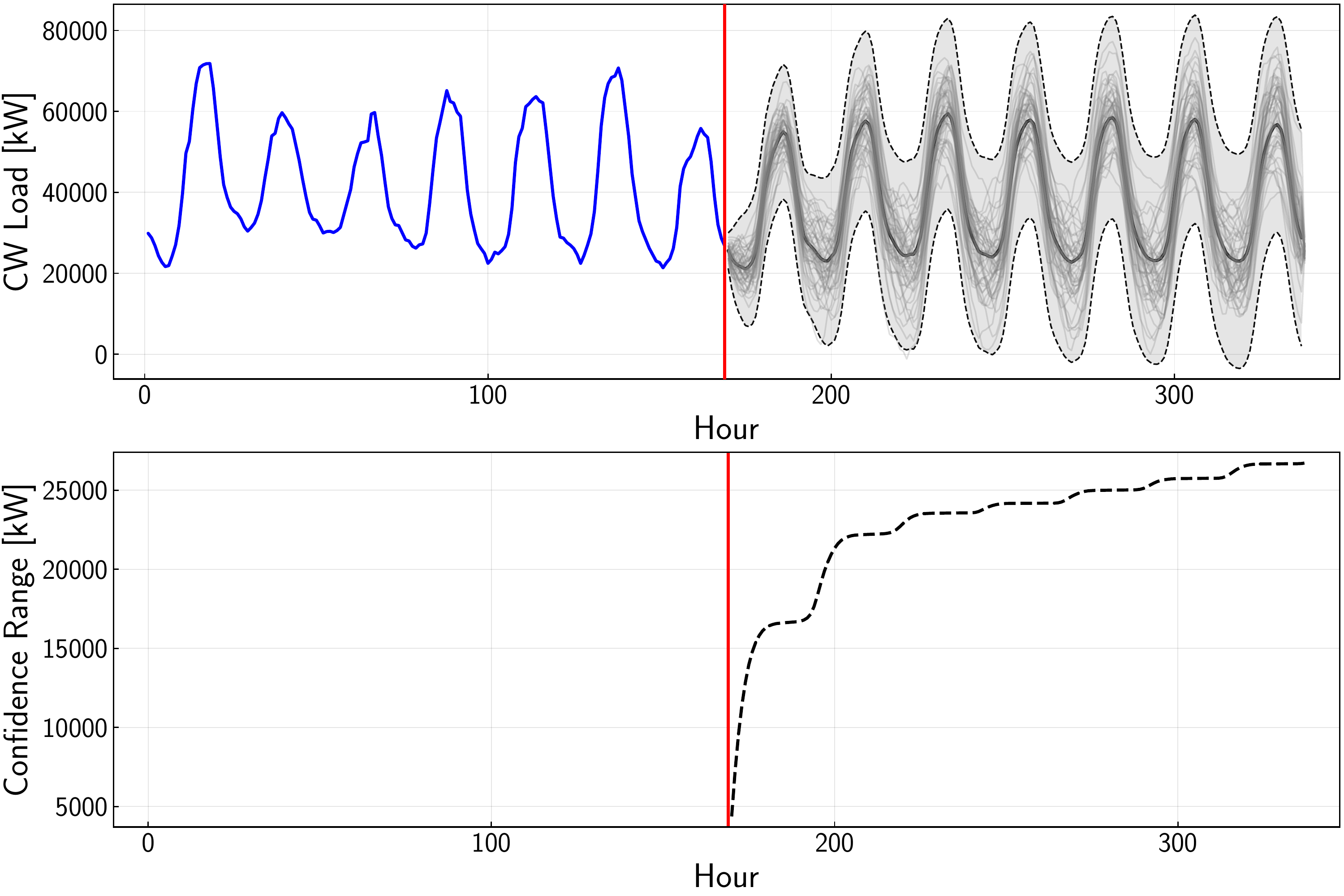}
\caption{A single instance of chilled water load forecast for 1-week with AR model. Vertical red line denotes the current time. In the top panel, dark bold curve represents mean forecast, the grey band denotes 99\% confidence interval and the light black curves within the band represent a few sample scenarios. The bottom panel shows the trajectory of the 99\% confidence range with prediction time.}
\label{fig:cw_fore}
\end{figure}

\begin{figure}[!htp]
\centering
\includegraphics[width=0.6\textwidth]{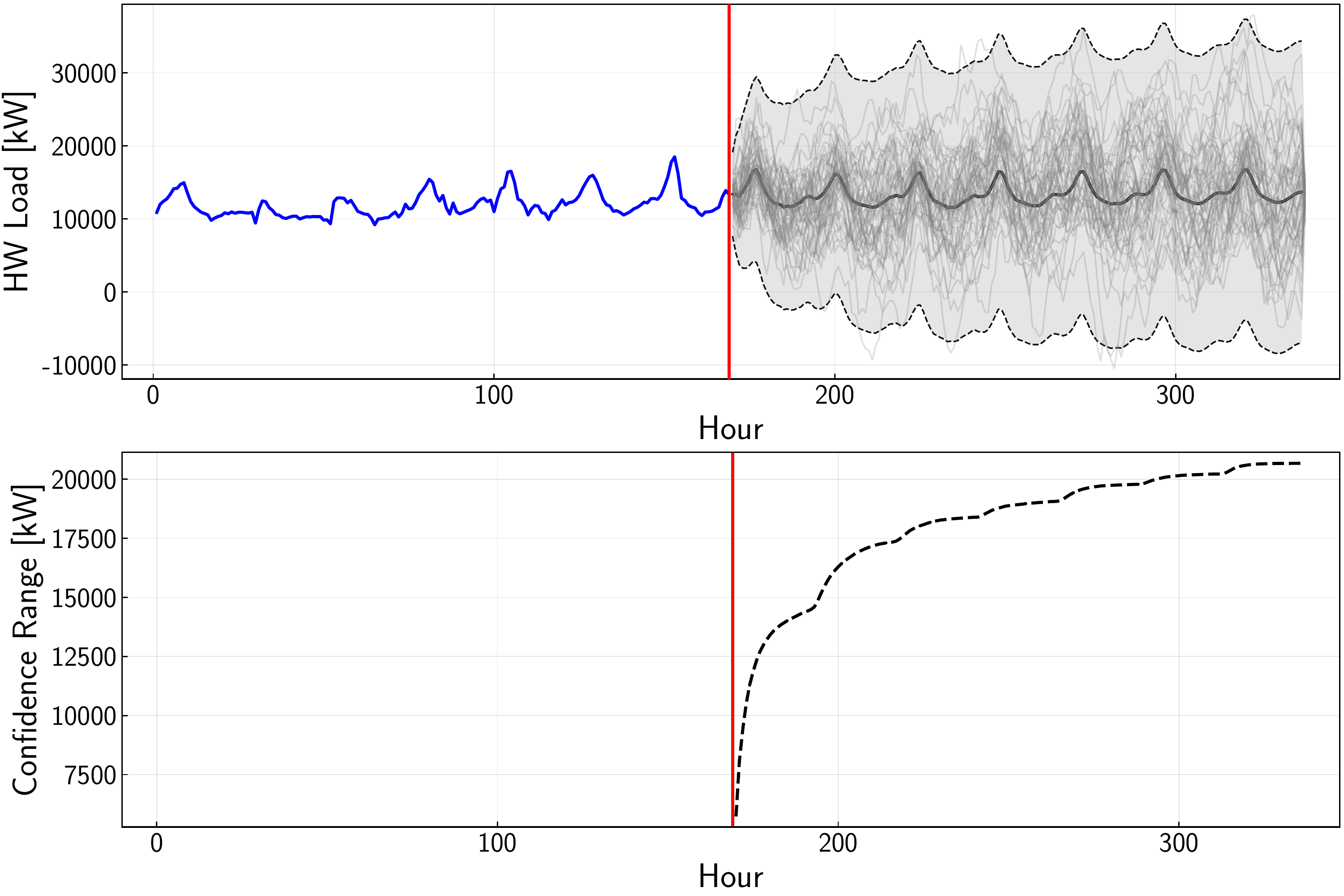}
\caption{A single instance of hot water load forecast for 1-week with AR model. Vertical red line denotes the current time. In the top panel, dark bold curve represents mean forecast, the grey band denotes 99\% confidence interval and the light black curves within the band represent a few sample scenarios. The bottom panel shows the trajectory of the 99\% confidence range with prediction time.}
\label{fig:hw_fore}
\end{figure}

\begin{figure}[!htp]
\centering
\includegraphics[width=0.6\textwidth]{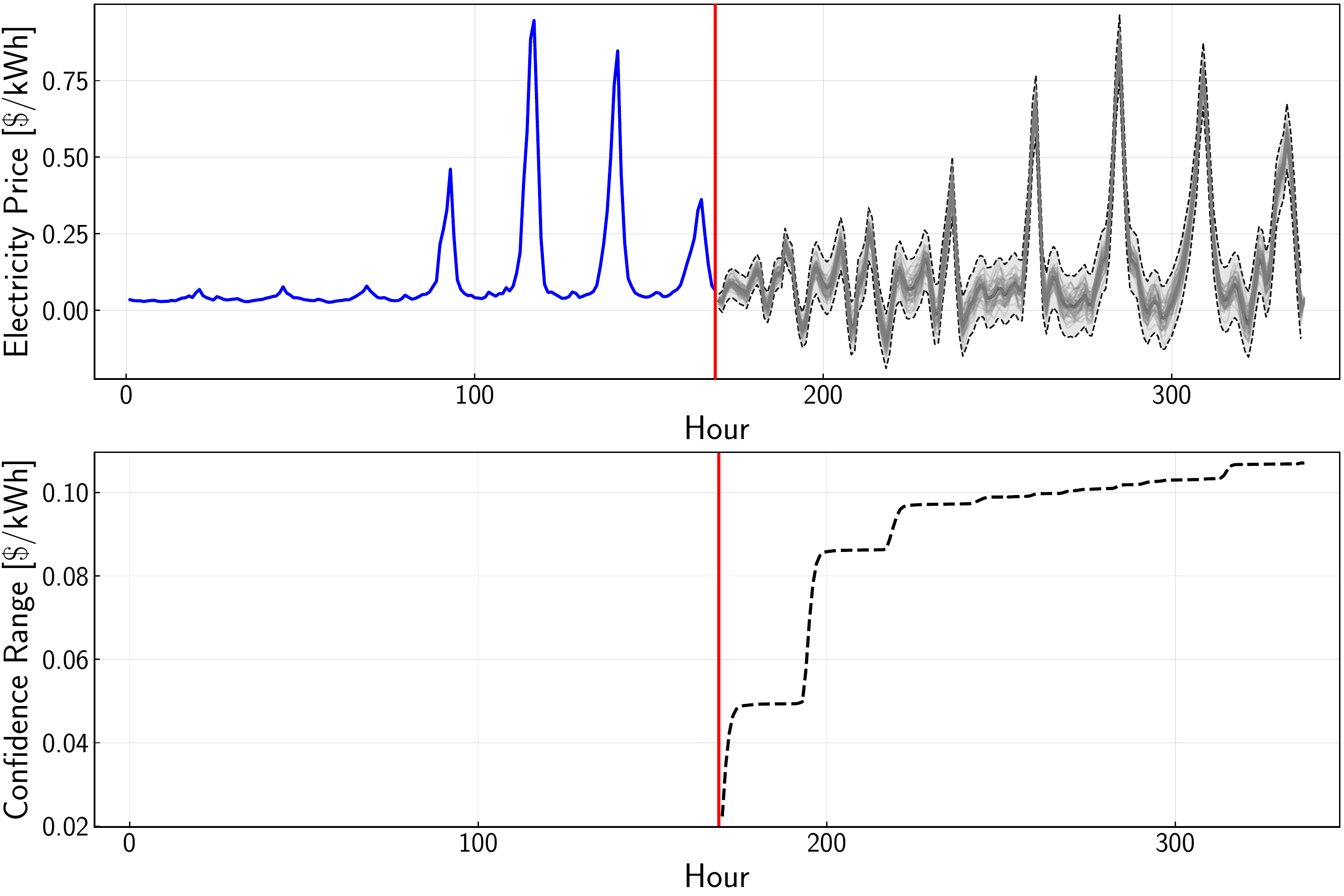}
\caption{A single instance of electricity price forecast for 1-week with AR model. Vertical red line denotes the current time. In the top panel, dark bold curve represents mean forecast, the grey band denotes 99\% confidence interval and the light black curves within the band represent a few sample scenarios. The bottom panel shows the trajectory of the 99\% confidence range with prediction time.}
\label{fig:price_fore}
\end{figure}

A single instance for the 1-week forecasts for electrical load, cold water load, hot water load, and electricity price is shown in Figures \ref{fig:load_fore}, \ref{fig:cw_fore}, \ref{fig:hw_fore} and \ref{fig:price_fore}. The mean forecasts are represented by the dark bold curves and the 99\% confidence intervals are shown in grey color. Sample scenarios are shown as light black curves. From these results we see that the AR models can capture the trends of the disturbances but that significant uncertainty exists. In particular, we notice that the magnitude of the confidence interval (the range) rises sharply within the first few hours. As we show next, this will be a major factor that drives constraint violations. 

\FloatBarrier

\FloatBarrier

\subsection{Closed-Loop Performance}

\begin{figure}[htp!]
\centering
\includegraphics[width=1\textwidth]{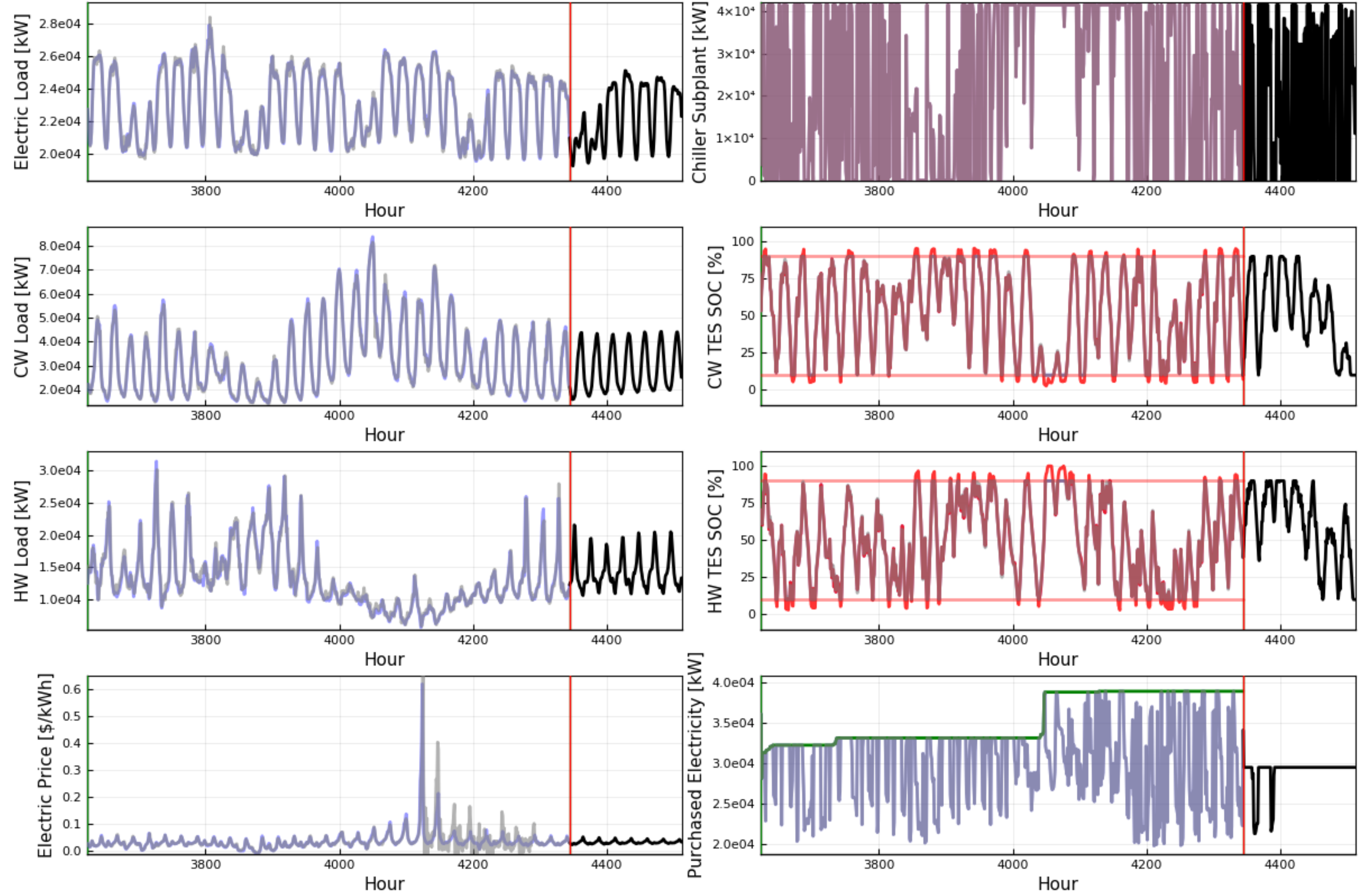}
\caption{Closed loop profile for deterministic MPC with $\beta = 0.1$. Black lines represent forecasts and model predictions. Blue lines represent actual realizations. For the control policies, red lines represent the actual implemented policy. For the state of charge, red horizontal lines represent the storage buffer. For the residual electrical load, the green line represents the running peak $R_t$.}
\label{fig:closed-loop-deterministic}
\end{figure}

\begin{figure}[htp!]
\centering
\includegraphics[width=1\textwidth]{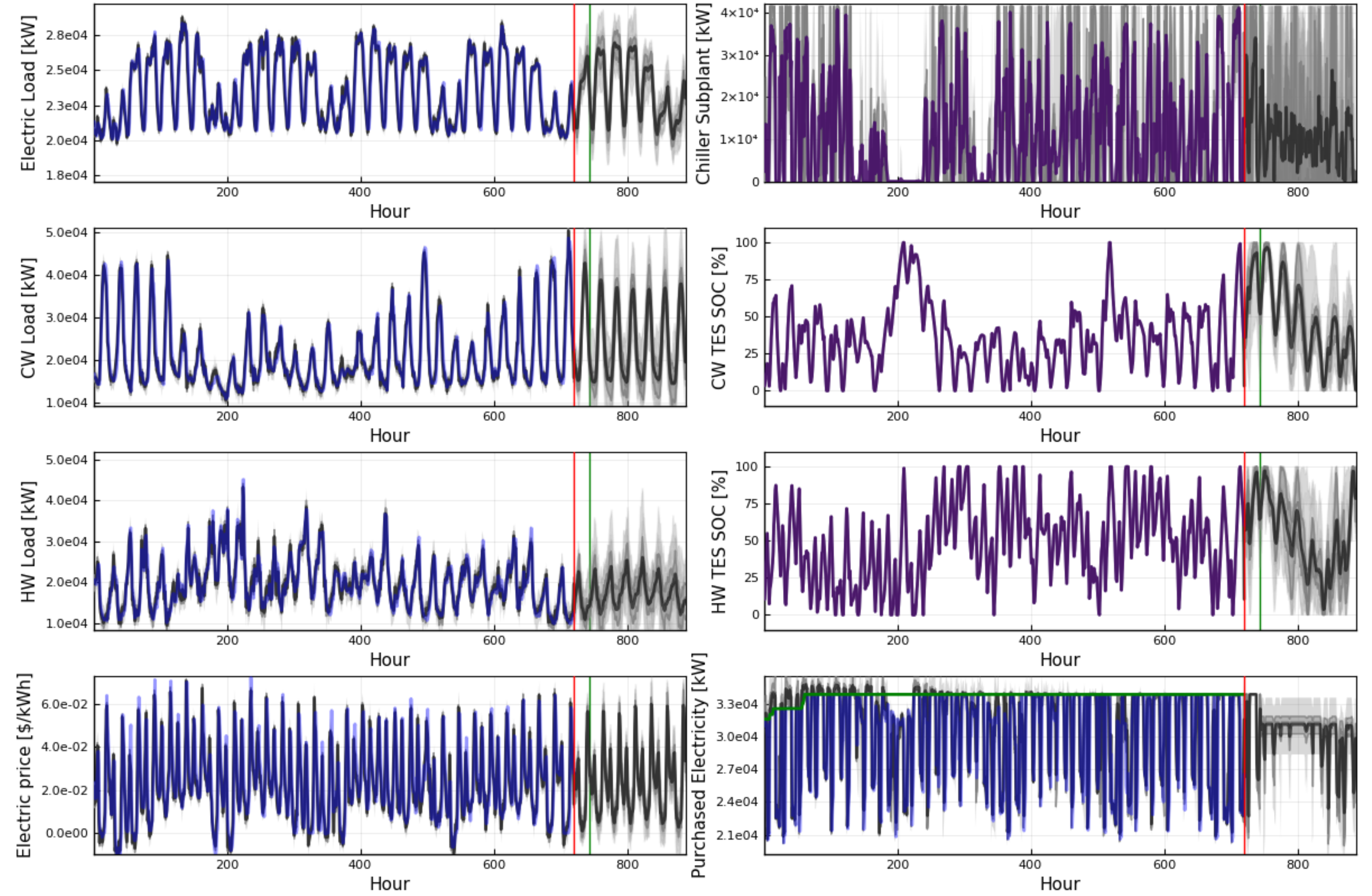}
\caption{Closed loop profile for stochastic MPC with $\beta = 0$. The grey regions represent uncertainty forecasts with black representing mean forecasts. Blue lines represent realized observations or committed policies. For the control policies, red lines represent the actual implemented policy and usually overlap with committed policy. For the residual load, the green line represents the peak observed residual electricity demand.}
\label{fig:closed-loop-stochastic}
\end{figure}

We compare the closed-loop policies for deterministic and stochastic MPC. In these results, a storage buffer of $\beta = 0$ is used for stochastic MPC and a buffer of 10\% ($\beta = 0.1$) is used for deterministic MPC. Figure \ref{fig:closed-loop-deterministic} provides a snapshot for a given validation scenario for deterministic MPC. Here, we note that the controller uses the storage buffer fairly frequently to counteract uncertainty in the disturbances. Figure \ref{fig:closed-loop-stochastic} shows a snapshot for stochastic MPC under a given validation scenario. Here, we note that the controller does not require an explicit buffer for storage, which results in a better utilization of the storage. An animation of the closed-loop performance of deterministic and stochastic MPC can be found at \url{https://github.com/zavalab/JuliaBox/tree/master/HVAC_Plant} (these help visualize the closed-loop dynamics). 

\FloatBarrier

\FloatBarrier

\subsection{Economic Performance and Constraint Violations}

In Table \ref{tab:utility_comparisons} we compare utility usage of the campus together with the HVAC central plant under the MPC implementations. We observe that the system would ideally (under perfect information) consume 264,800 MWh of electricity, 197,839 million gallons of water, and 29,067 MWh of natural gas in the year. These numbers represent the best possible performance \mynote{under a finite horizon MPC} and highlight the high resource use of these systems. The system under a deterministic MPC implementation has a higher consumption for each of the utilities. This is expected because deterministic MPC faces forecast errors. The stochastic MPC implementation reduces natural gas use of the deterministic counterpart by a significant amount. Specifically, under stochastic MPC, we observe reductions in natural gas use of 8.57\% (2,840 MWh).  These results highlight that stochastic MPC can reduce the resource footprint by better handling of disturbances. Stochastic MPC achieves modest reductions in electricity and water consumption of 0.07\% (187 MWh) and 0.96\% (1,927 million gallons), respectively. We will see, however, that significant cost reductions are achieved (due to the time-varying nature of electricity prices). 

\begin{table*}[htp!]\centering
\caption{Utility usage analysis for different MPC implementations.}
\resizebox{\columnwidth}{!}{%
\begin{tabular}{|C{4cm}|c|c|c|c|}
\hline
\multirow{2}{*}{Expected Utility} & \multicolumn{4}{c|}{\textbf{Campus With Central Plant}}\\ \cline{2-5}
{Usage}  & {\begin{tabular}{@{}c@{}}\text{Perfect} \\ \text{Information}\end{tabular}} & {\begin{tabular}{@{}c@{}}\text{Stochastic} \\{$\beta = 0$}\end{tabular}} & {\begin{tabular}{@{}c@{}}\text{Deterministic} \\{$\beta = 0.1$}\end{tabular}}
& {\begin{tabular}{@{}c@{}}\text{Savings by} \\ \text{Stochastic MPC}\end{tabular}} \\
\hline
{\begin{tabular}{@{}c@{}}\textbf{Electricity} \\{(MWh/year)}\end{tabular}}  & {\begin{tabular}{@{}c@{}} {264,800} \\ \textbf{(-0.2\%)}\end{tabular}} & {\begin{tabular}{@{}c@{}} {265,137} \\ \textbf{(-0.07\%)}\end{tabular}} & {\begin{tabular}{@{}c@{}} {265,324} \\ \textbf{(Base)}\end{tabular}} & {187} \\
\hline
{\begin{tabular}{@{}c@{}}\textbf{Water} \\ {(MMgal/year)}\end{tabular}}   & {\begin{tabular}{@{}c@{}} {{197,839}} \\ \textbf{(-1.8\%)}\end{tabular}} &  {\begin{tabular}{@{}c@{}} {199,593} \\ \textbf{(-0.96\%)}\end{tabular}} & {\begin{tabular}{@{}c@{}} {201,520} \\ \textbf{(Base)}\end{tabular}} & {1,927} \\
\hline
{\begin{tabular}{@{}c@{}}\textbf{Natural Gas} \\ {(MWh/year)}\end{tabular}}  & {\begin{tabular}{@{}c@{}} {29,067} \\ \textbf{(-12.3\%)}\end{tabular}} &  {\begin{tabular}{@{}c@{}} {30,315} \\ \textbf{(-8.57\%)}\end{tabular}} & {\begin{tabular}{@{}c@{}} {33,155} \\ \textbf{(Base)}\end{tabular}} & {2,840} \\
\hline
\end{tabular} \label{tab:utility_comparisons}
}
\end{table*}
\FloatBarrier

The economic performance of the MPC implementations is summarized in Table \ref{tab:withcentralplant_comparisons}. Here, we present the expected total cost and we disaggregate the this cost in its different components. The expected total cost of the campus (without the central plant) is 11,815,567 \$/year and the expected total cost for the campus with the central plant (operated with perfect information MPC) is 16,047,162 \$/year. This indicates that the operation of the HVAC central plant alone costs 4,231,595 \$/year (this shows the large costs associated with the plant). We factor out the cost of the central plant from the total cost because the MPC controllers can only help reduce the central plant costs (the campus costs are exogenous). 

From Table \ref{tab:withcentralplant_comparisons} we also see that the expected cost of the central plant is improved by 7.52\% by using stochastic MPC (relative to deterministic MPC). The associated cost savings total 349,000 \$/year. These savings represent 75\% of the possible improvement over deterministic MPC (obtained with perfect information MPC). 

Table \ref{tab:withoutcentralplant_comparisons} disaggregates the costs of the central plant.  Here, we observe that stochastic MPC  achieves an improvement of 29.8\% in the demand charge cost over deterministic MPC and achieves improvements in electricity cost of 6.88\%, in natural gas cost of 8.57\%, and in water cost of 0.96\%. We thus see that, even if reductions in electricity use are moderate, reductions in cost are significant. We also note that reductions in natural gas use and cost are both 8.57\% (because the price of gas is not varying over time). The dramatic reduction in demand charge costs indicates that disturbance uncertainty has a strong effect on peak electricity load (we recall that disturbance uncertainty rises sharply over time).

\begin{table*}[htp!]\centering
\caption{Economic performance analysis for different MPC implementations (central plant only).}
\resizebox{\columnwidth}{!}{%
\begin{tabular}{|C{4cm}|c|c|c|c|c|}
\hline
\multirow{2}{*}{} & \multirow{2}{*}{\textbf{Campus Only}} & \multicolumn{4}{c|}{\textbf{Central Plant Only}}\\ \cline{3-6}
{} & {(No Central Plant)} & {\begin{tabular}{@{}c@{}}\text{Perfect} \\ \text{Information}\end{tabular}} & {\begin{tabular}{@{}c@{}}\text{Stochastic} \\{$\beta = 0$}\end{tabular}} & {\begin{tabular}{@{}c@{}}\text{Deterministic} \\{$\beta = 0.1$}\end{tabular}}
& {\begin{tabular}{@{}c@{}}\text{Value of} \\ \text{Stochastic MPC}\end{tabular}} \\
\hline
{\begin{tabular}{@{}c@{}}\textbf{Expected Electricity}\\ \textbf{Cost} \\{(\$/year)}\end{tabular}} & {9,764,251} & {1,605,815} &  {1,658,335} & {1,780,771} & {\begin{tabular}{@{}c@{}}{122,436}\\ \textbf{(6.88\%)} \end{tabular}} \\
\hline
{\begin{tabular}{@{}c@{}}\textbf{Expected Water Cost} \\ {(\$/year)}\end{tabular}} & {-} & {1,780,554} &  {1,796,342} & {1,813,681} & {\begin{tabular}{@{}c@{}}{17,339}\\ \textbf{(0.96\%)} \end{tabular}} \\
\hline
{\begin{tabular}{@{}c@{}}\textbf{Expected Natural Gas} \\ \textbf{Cost} \\ {(\$/year)}\end{tabular}} & {-} & {523,210} &  {545,679} & {596,801} &  {\begin{tabular}{@{}c@{}}{51,122}\\ \textbf{(8.57\%)} \end{tabular}}  \\
\hline
{\begin{tabular}{@{}c@{}}\textbf{Expected Demand} \\ \textbf{Charge} \\ {(\$/year)}\end{tabular}} & {2,051,316} & {222,015} &  {315,637} & {449,206} & {\begin{tabular}{@{}c@{}}{133,569}\\ \textbf{(29.8\%)} \end{tabular}}  \\
\hline
{\begin{tabular}{@{}c@{}}\textbf{Cost of} \\ \textbf{Central Plant} \\ {(\$/year)}\end{tabular}} & {-} & {4,231,595} &  {4,291,497} & {4,640,460} & {\begin{tabular}{@{}c@{}}{348,963}\\ \textbf{(7.52\%)} \end{tabular}}  \\
\hline
{\begin{tabular}{@{}c@{}}\textbf{Savings in} \\ \textbf{Central Plant Cost}\end{tabular}} & {-} & {9.66\%} &  {7.52\%} & \textbf{Base} & {7.52\%} \\
\hline
\end{tabular} \label{tab:withoutcentralplant_comparisons}
}
\end{table*}
\FloatBarrier

\begin{figure}[htp!]
\centering
\includegraphics[width=\textwidth]{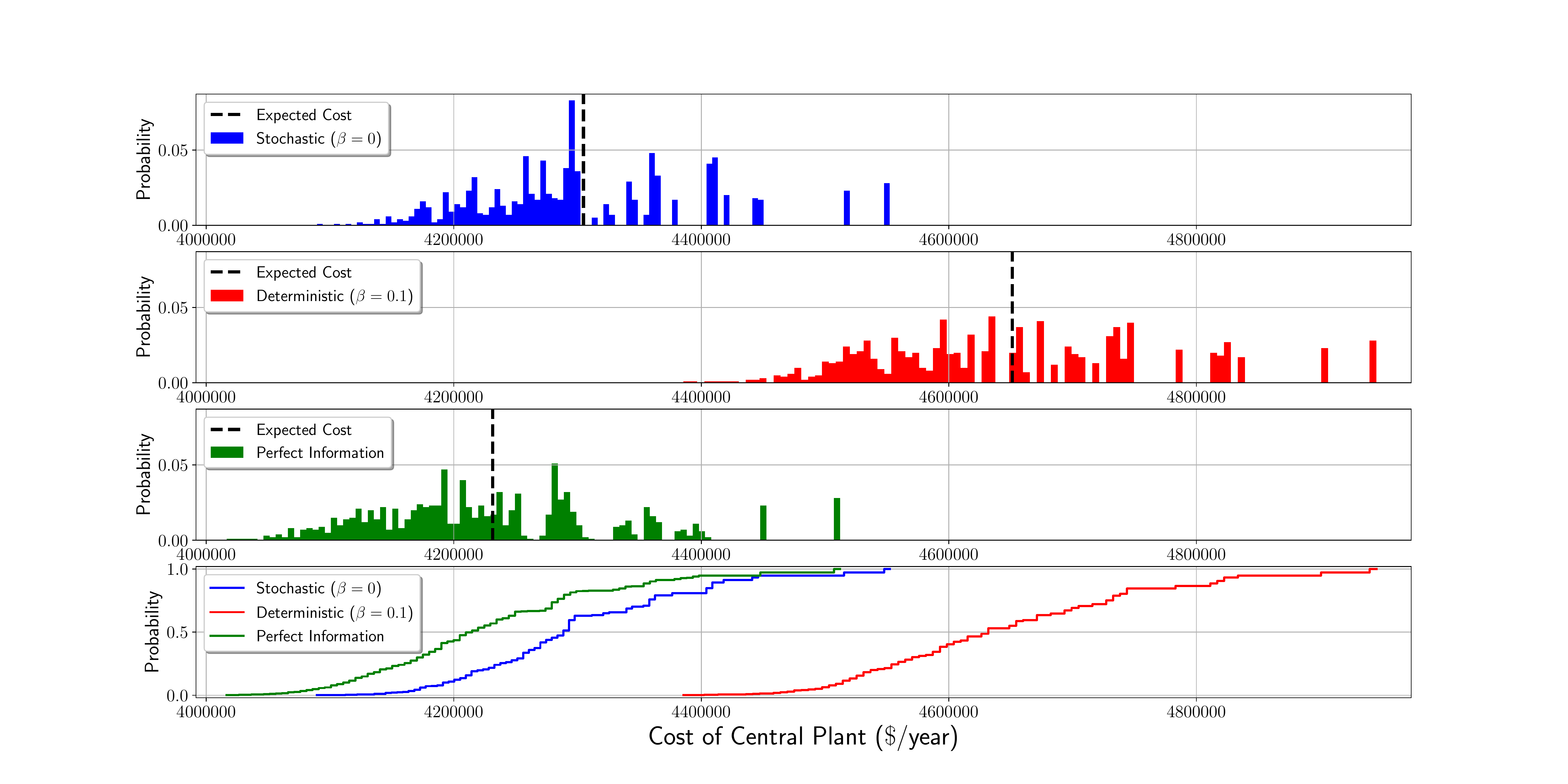}
\caption{Probability and cumulative distributions for total cost of central plant.}
\label{fig:costs}
\end{figure}

\begin{figure}[htp!]
\centering
\includegraphics[width=\textwidth]{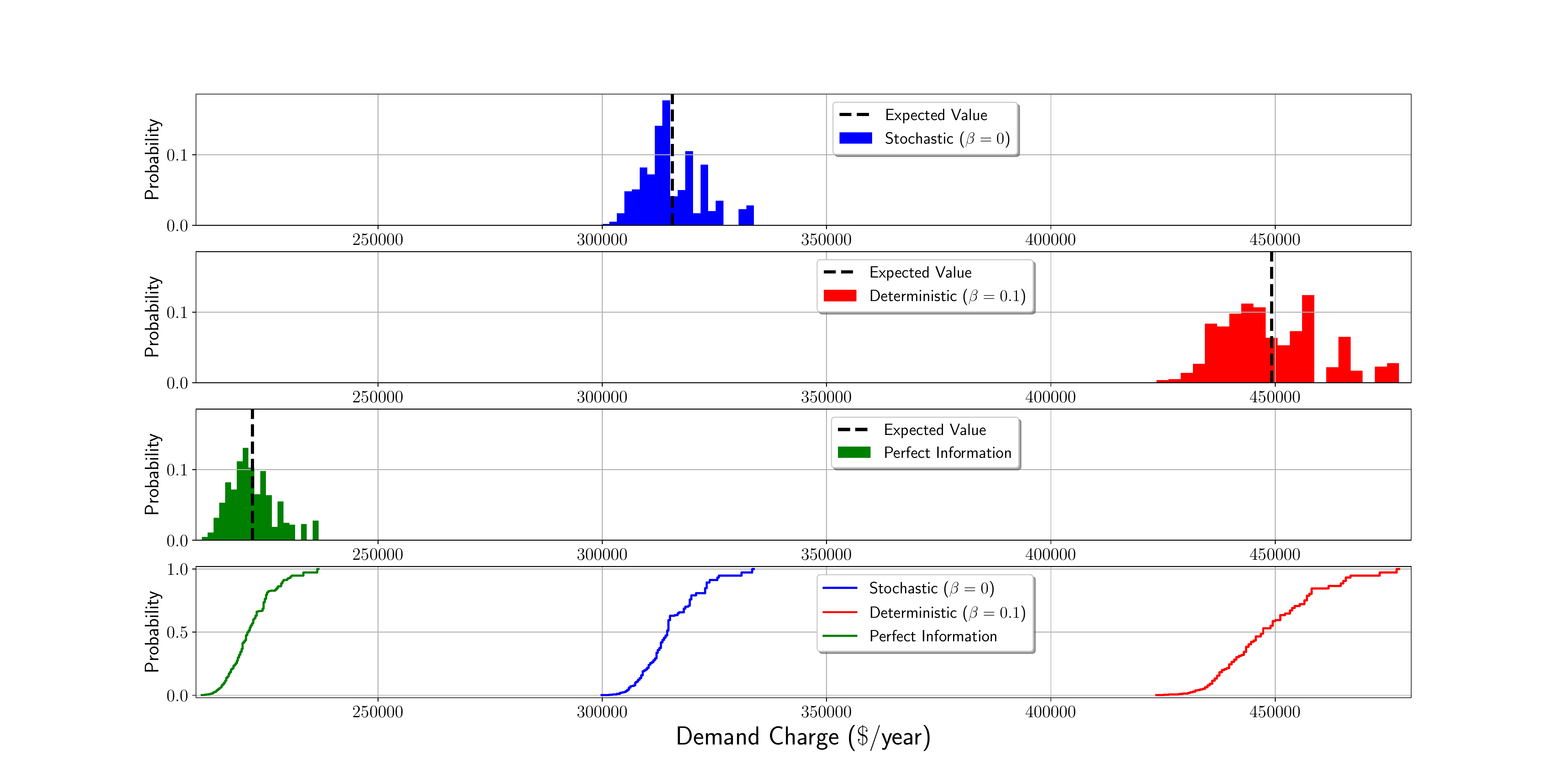}
\caption{Probability and cumulative distributions for demand charges.}
\label{fig:peaks}
\end{figure}

The probability distribution and cumulative distribution for the costs of the HVAC central plant and for the peak demand charges are shown in Figures \ref{fig:costs} and \ref{fig:peaks}. From the cumulative distributions we observe that the probability of taking a smaller cost and a smaller demand charge for the central plant is higher under stochastic MPC compared to deterministic MPC. Similarly, the probability of obtaining a higher cost and a higher demand charge for the central plant is much higher under deterministic MPC than under stochastic MPC. This clearly illustrates that the performance of stochastic MPC consistently dominates that of deterministic MPC.

\begin{figure}[htp!]
\centering
\includegraphics[width=0.7\textwidth]{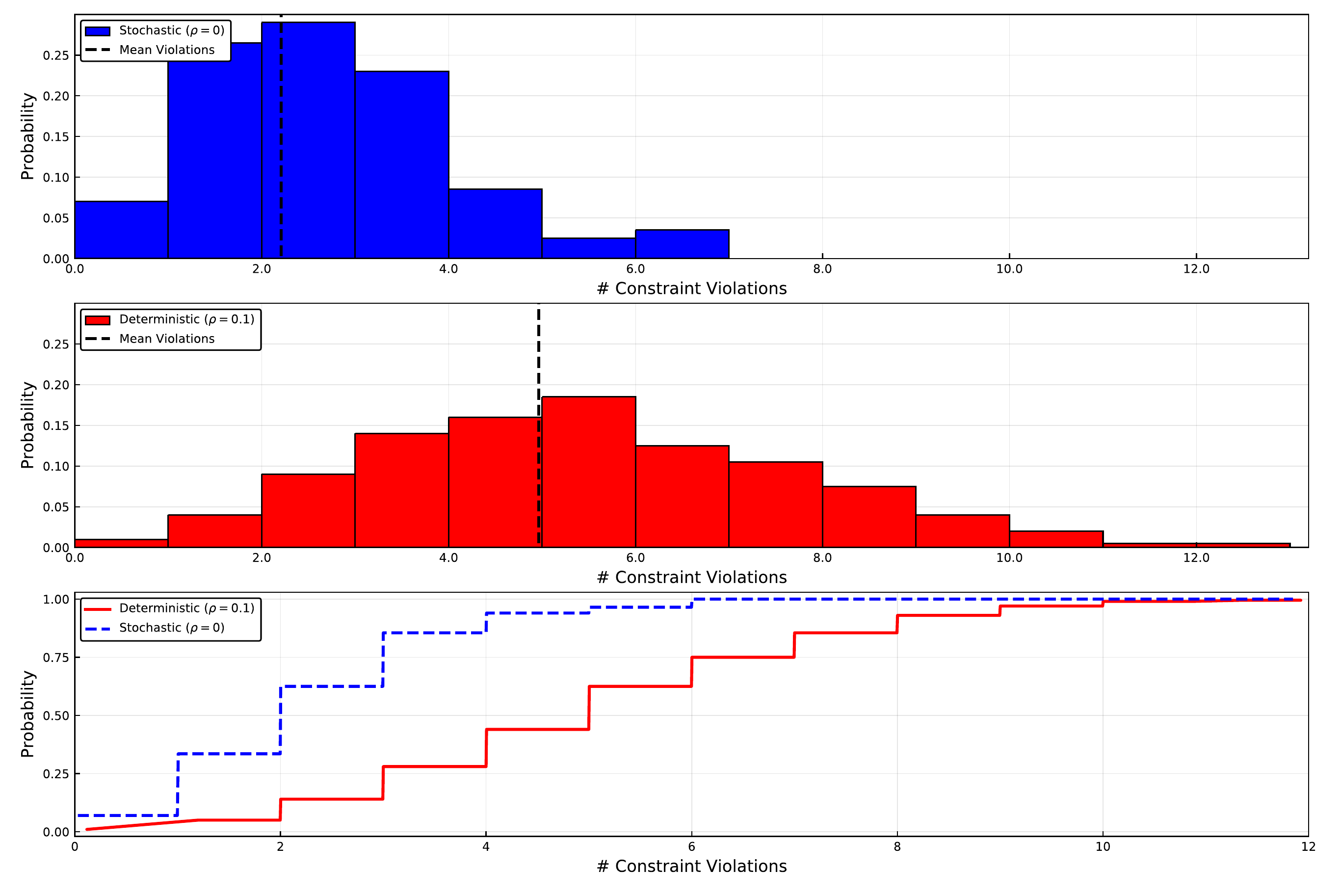}
\caption{Probability and cumulative distributions for constraint violations (per 100 hours).}
\label{fig:violations}
\end{figure}

We then evaluated constraint violations (in terms of storage overflow or drying up) obtained with deterministic and stochastic MPC for all the validation scenarios. We recall that the MPC controllers take no action when the control policy becomes infeasible. This can result in either overflow or drying up of the chilled water and hot water storage tanks, and therefore can lead to loss of energy. Figure \ref{fig:violations} shows how often infeasibility occurs per 100 hours over the year for the 200 validation scenarios. The results highlight that stochastic MPC without a storage buffer is more reliable at maintaining a feasible operation than deterministic MPC with a 10\% buffer. 

\FloatBarrier

\begin{figure}[htp!]
\centering
\includegraphics[width=0.7\textwidth]{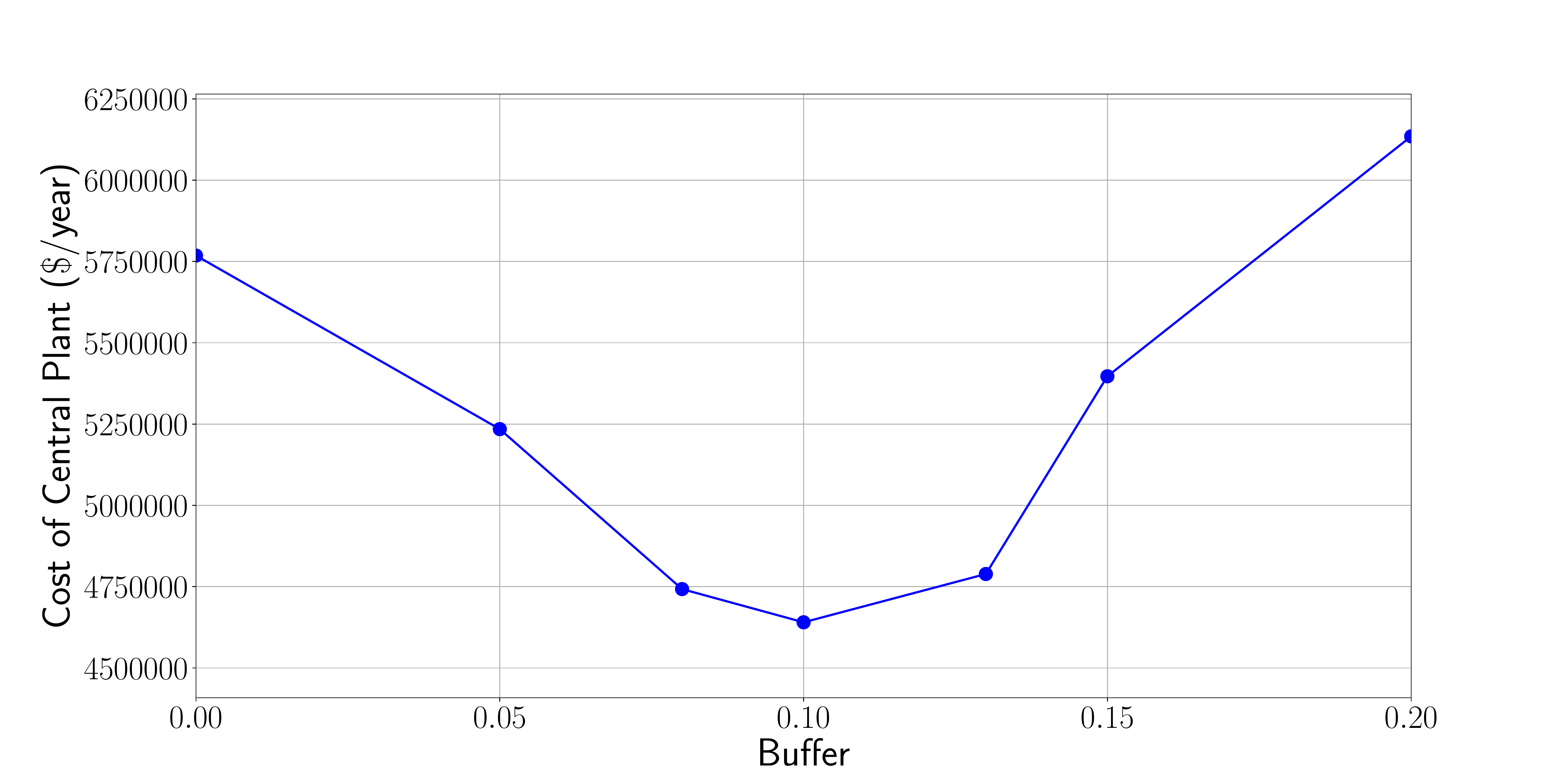}
\caption{Economic performance (cost of central plant) of deterministic MPC with varying buffers.}
\label{fig:buffer_sensitivity_CP}
\end{figure}

In the above case study, we chose a buffer of 10\% for deterministic MPC as this was the best possible buffer value found.   To see this, we examined the effect of using different buffers on economic performance. Table \ref{tab:sensitivity_buffer_deterministic} compares results for buffers of 0-20\% for deterministic MPC. As can be seen, deterministic MPC performs poorly without any buffer and results in a very high cost (compared to performance obtained with a buffer of 10\%). This is because a deterministic MPC implementation with 0\% storage buffer tries to utilize the full capacity of the TES in each optimization step but, after the actual realization of the loads is observed, it has to adjust the control actions because of frequent constraint violations which in turn leads to higher demand charges. Increasing the buffer initially leads to lower total cost for deterministic MPC because the controller is able to avoid infeasibility issues but eventually becomes detrimental because the fraction of the storage capacity available decreases (thus limiting flexibility). This inherent trade-off is shown in Table  \ref{tab:sensitivity_buffer_deterministic} and Figures \ref{fig:buffer_sensitivity_CP}. Here, it is clear that a buffer of 10\% achieves the best cost. We note, however, that even this best tuned cannot achieve the performance of stochastic MPC. In other words, stochastic MPC consistently dominates deterministic MPC. 

\begin{table*}[htp!]\centering
\caption{Expected costs for deterministic MPC with varying buffers.}
\resizebox{\columnwidth}{!}{%
\begin{tabular}{|c|c|c|c|c|c|c|c|}
\hline
\textbf{Item} & {\begin{tabular}{@{}c@{}}\textbf{Deterministic}\\{($\rho = 0$)}\end{tabular}} & {\begin{tabular}{@{}c@{}}\textbf{Deterministic}\\{($\beta = 0.05$)}\end{tabular}} & {\begin{tabular}{@{}c@{}}\textbf{Deterministic}\\{($\beta = 0.08$)}\end{tabular}} & {\begin{tabular}{@{}c@{}}\textbf{Deterministic}\\{($\beta = 0.1$)}\end{tabular}} & {\begin{tabular}{@{}c@{}c@{}}\textbf{Deterministic}\\{($\beta = 0.13$)}\end{tabular}} & {\begin{tabular}{@{}c@{}}\textbf{Deterministic}\\{($\beta = 0.15$)}\end{tabular}} & {\begin{tabular}{@{}c@{}}\textbf{Deterministic}\\{($\beta = 0.2$)}\end{tabular}}\\
\hline
 {\begin{tabular}{@{}c@{}}\textbf{Total Cost}\\{(\$/year)}\end{tabular}} & 17,583,658  & 17,050,005 & 16,558,117 & 16,456,027 & 16,604,521 & 17,212,421 & 17,950,255\\
\hline
{\begin{tabular}{@{}c@{}}\textbf{Cost of Central Plant}\\{(\$/year)}\end{tabular}} & 5,768,091  & 5,234,438 & 4,742,550 & 4,640,461  & 4,788,954 & 5,396,854 & 6,134,688 \\
 \hline
{\begin{tabular}{@{}c@{}}\textbf{Improvement in} \\ \textbf{Cost of Central Plant}\end{tabular}} & -24.3\% & -12.8\% & -2.20\% & \textbf{Base} & -3.21\% & -16.3\% & -32.2\%\\
\hline
\end{tabular} \label{tab:sensitivity_buffer_deterministic}
}
\end{table*}

\FloatBarrier

\section{Conclusions and Future Work}

We presented a stochastic MPC framework for a HVAC central plant of a typical university campus. We use the framework to rigorously assess of the benefits of stochastic MPC over deterministic MPC in terms of economic performance and constraint violations. Our framework uses real historical data to conduct forecasting and uncertainty quantification tasks. Our results show that stochastic MPC reduces overall energy consumption and cost by better handling of storage and better integration of electricity, natural gas, and water. Specifically, we found that stochastic MPC can reduce the natural gas cost by 8.57\%, electricity cost by 6.89\%, and peak demand charges by 29.8\% (relative to deterministic MPC). We find that deterministic MPC leads to frequent constraint violations in storage capacity (causing overflow or drying up) \mynote{of the hot and chiller water tanks}. Stochastic MPC is able to avoid constraint violations because it anticipates the uncertainty by explicitly incorporating it in the model.  As part of future work, we will seek to implement parallel decomposition and simulations in order to accelerate simulations. \mynote{The use of affine decision rules and multi-stage formulations are also interesting future directions.}

\section{Acknowledgements}
We acknowledge funding from the National Science Foundation under award NSF-EECS-1609183. VZ declares a financial interest in Johnson Controls International, a for-profit company that develops control technologies. 

\section*{APPENDIX} \label{sec:appendix}
\appendix
\section{Additional Results}
\begin{table*}[htp!]\centering
\caption{Economic performance analysis for different MPC implementations (campus with central plant).}
\resizebox{\columnwidth}{!}{%
\begin{tabular}{|C{4cm}|c|c|c|c|c|}
\hline
\multirow{2}{*}{} & \multirow{2}{*}{\textbf{Campus Only}} & \multicolumn{4}{c|}{\textbf{Campus With Central Plant}}\\ \cline{3-6}
{} & {(No Central Plant)} & {\begin{tabular}{@{}c@{}}\text{Perfect} \\ \text{Information}\end{tabular}} & {\begin{tabular}{@{}c@{}}\text{Stochastic} \\{$\beta = 0$}\end{tabular}} & {\begin{tabular}{@{}c@{}}\text{Deterministic} \\{$\beta = 0.1$}\end{tabular}}
& {\begin{tabular}{@{}c@{}}\text{Value of} \\ \text{Stochastic MPC}\end{tabular}} \\
\hline
{\begin{tabular}{@{}c@{}}\textbf{Expected Electricity}\\ \textbf{Cost} \\{(\$/year)}\end{tabular}} & {9,764,251} & {11,370,066} &  {11,422,586} & {11,545,022} & 122,436 \\
\hline
{\begin{tabular}{@{}c@{}}\textbf{Expected Water Cost} \\ {(\$/year)}\end{tabular}} & {-} & {1,780,554} &  {1,796,342} & {1,813,681} & {17,339} \\
\hline
{\begin{tabular}{@{}c@{}}\textbf{Expected Natural Gas} \\ \textbf{Cost} \\ {(\$/year)}\end{tabular}} & {-} & {523,210} &  {545,679} & {596,801} & {51,122} \\
\hline
{\begin{tabular}{@{}c@{}}\textbf{Expected Demand} \\ \textbf{Charge} \\ {(\$/year)}\end{tabular}} & {2,051,316} & {2,273,331} &  {2,342,457} & {2,500,522} & {133,569} \\
\hline
{\begin{tabular}{@{}c@{}}\textbf{Expected Total} \\ \textbf{Cost} \\ {(\$/year)}\end{tabular}} & {11,815,567} & {16,047,162} &  {16,107,064} & {16,456,027} & {348,963} \\
\hline
{\begin{tabular}{@{}c@{}}\textbf{Cost of} \\ \textbf{Central Plant} \\ {(\$/year)}\end{tabular}} & {-} & {4,231,595} &  {4,291,497} & {4,640,460} & {348,963} \\
\hline
{\begin{tabular}{@{}c@{}}\textbf{Savings in} \\ \textbf{Central Plant Cost}\end{tabular}} & {-} & {9.66\%} &  {7.52\%} & \textbf{Base} & {7.52\%} \\
\hline
\end{tabular} \label{tab:withcentralplant_comparisons}
}
\end{table*}
\FloatBarrier
\section{Nomenclature} \label{subsec:nomen}
\noindent \textbf{Sets and indices:}
\begin{itemize}
\item[\textbullet] $\mathcal{T} :=  \{t,t+1,t+2,\dots,t+N-1\}$: Prediction horizon time set, where $N$ is the prediction horizon length, $t$ is the current time.
\item[\textbullet] $\mathcal{Y} :=  \{1,2,\dots,Y\}$: Planning horizon time set, where $Y$ is the planning (or simulation) horizon length.
\item[\textbullet] $\mathcal{H} :=  \{t-H,t-H+1,\dots,t\}$: Historical horizon set, with $H$ being the length.
\item[\textbullet] $\mathcal{T}_M :=  \{t,t+1,t+2,\dots,t+N-1\}$: Set of ending time indices for all months. 
\item[\textbullet] $\bar{\Xi}$: Set of scenarios in stochastic MPC formulation.
\item[\textbullet] $\tilde{\Xi} $: Set of scenarios used for validation.
\item[\textbullet] $t$: Time instant index.
\item[\textbullet] $t_m$: Ending time index (hour) of month $m$.
\item[\textbullet] $\xi$: Realization index.
\end{itemize}

\noindent \textbf{Model Parameters and Data:}
\begin{itemize}
\item[\textbullet] $\pi^{e}_t \textrm{ and } \pi^{e}_t(\xi) \in \mathbb{R}_+$: Electricity price [\$/kWh] over the time interval $[t,(t+1)]$.
\item[\textbullet] $\hat{\pi}^{e}_t \in \mathbb{R}_+$: Forecast electricity price [\$/kWh].
\item[\textbullet] $\pi^{w}_t \in \mathbb{R}_+$: Price of water [\$/gal] over the time interval $[t,(t+1)]$.
\item[\textbullet] $\pi^{ng}_t \in \mathbb{R}_+$: Price of natural gas [\$/kWh] over the time interval $[t,(t+1)]$.
\item[\textbullet] $\pi^D \in \mathbb{R}_+$: Rate of demand charge [\$/kW].
\item[\textbullet] $\alpha^{e}_{cs} \in \mathbb{R}_+$: kW of electricity used by chiller subplant per kW chilled water produced [-].
\item[\textbullet] $\alpha^{e}_{hrc} \in \mathbb{R}_+$: kW of electricity used by HR chiller subplant per kW chilled water produced [-].
\item[\textbullet] $\alpha^{e}_{hwg} \in \mathbb{R}_+$: kW of electricity used by hot water generator per kW hot water produced [-].
\item[\textbullet] $\alpha^{e}_{ct} \in \mathbb{R}_+$: kW of electricity used by cooling towers per kW condenser water input [-].
\item[\textbullet] $\alpha^{w}_{ct} \in \mathbb{R}_+$: Gallons of water used by cooling towers per kW condenser water input [-].
\item[\textbullet] $\alpha^{ng}_{hwg} \in \mathbb{R}_+$: kW of natural gas used by hot water generator per kW hot water produced [-].
\item[\textbullet] $\alpha^{cond}_{cs} \in \mathbb{R}_+$: kW of condenser water produced by chiller subplant per kW chilled water produced [-].
\item[\textbullet] $\alpha^{h}_{hrc} \in \mathbb{R}_+$: kW of hot water produced by HR chiller subplant per kW chilled water produced [-].
\item[\textbullet] $\rho^{cw} \in \mathbb{R}_+$: Penalty for unmet chilled water load [\$/kWh].
\item[\textbullet] $\rho^{hw} \in \mathbb{R}_+$: Penalty for unmet hot water load [\$/kWh].
\item[\textbullet] $\sigma_t := \min\{(M-t)/N,1\}$: Discounting factor for the monthly demand charge price [-].
\item[\textbullet] $\beta \in [0,1]$: Storage buffer for the chilled water and hot water TES [-].
\item[\textbullet] $L^{e}_{t} \textrm{ and } L^{e}_t(\xi) \in \mathbb{R}_+$: Electrical load of campus [kW] over the time interval $[t,(t+1)]$.
\item[\textbullet] $\hat{L}^{e}_t \in \mathbb{R}_+$: Forecast electrical load [kW].
\item[\textbullet] $L^{cw}_{t} \textrm{ and } L^{cw}_t(\xi) \in \mathbb{R}_+$: Chilled water load [kW] over the time interval $[t,(t+1)]$.
\item[\textbullet] $\hat{L}^{cw}_t \in \mathbb{R}_+$: Forecast chilled water load [kW].
\item[\textbullet] $L^{hw}_{t} \textrm{ and } L^{hw}_t(\xi) \in \mathbb{R}_+$: Hot water load [kW] over the time interval $[t,(t+1)]$.
\item[\textbullet] $\hat{L}^{cw}_t \in \mathbb{R}_+$: Forecast hot water load [kW].
\item[\textbullet] $\overline{E}_{cw}\in \mathbb{R}_+$: Energy storage capacity of chilled water energy storage [kWh].
\item[\textbullet] $\overline{E}_{hw}\in \mathbb{R}_+$: Energy storage capacity of hot water energy storage [kWh].
\item[\textbullet] $\overline{P}_{cs}\in \mathbb{R}_+$: Maximum load of chiller subplant [kW].
\item[\textbullet] $\overline{P}_{hrc}\in \mathbb{R}_+$: Maximum load of heat recovery (HR) chiller subplant [kW].
\item[\textbullet] $\overline{P}_{hwg}\in \mathbb{R}_+$: Maximum load of hot water generator [kW].
\item[\textbullet] $\overline{P}_{ct}\in \mathbb{R}_+$: Maximum load of cooling towers [kW].
\item[\textbullet] $\overline{P}_{cw}\in \mathbb{R}_+$: Maximum discharging rate of chilled water energy storage [kW].
\item[\textbullet] $\overline{P}_{hw}\in \mathbb{R}_+$: Maximum discharging rate of hot water energy storage [kW].
\item[\textbullet] $R_t\in \mathbb{R}_+$: Peak electrical load observed until time $t\in\mathcal{T}_M$ [kW].
\item[\textbullet] $v_{cw,t}$: Random variable used to update chilled water storage from predicted value to simulated actual value.
\item[\textbullet] $v_{hw,t}$: Random variable used to update hot water storage from predicted value to simulated actual value.
\item[\textbullet] $\sigma^2_{cw,err,t+1}$: Variance of chilled water load prediction error for $t+1$.
\item[\textbullet] $\sigma^2_{hw,err,t+1}$: Variance of hot water load prediction error for $t+1$.
\item[\textbullet] $\sigma^2_{cw,int}$: Variance of integrated chilled water load for 1-hour periods.
\item[\textbullet] $\sigma^2_{hw,int}$: Variance of integrated hot water load for 1-hour periods.
\end{itemize}

\noindent \textbf{Controls:}
\begin{itemize}
\item[\textbullet] $P_{cs,t} \in \mathbb{R}_+$: Amount of chilled water produced by chiller subplant [kW] over the time interval $[t,(t+1)]$. 
\item[\textbullet] $P_{hrc,t} \in \mathbb{R}_+$: Amount of chilled water produced by heat recovery (HR) chiller subplant [kW] over the time interval $[t,(t+1)]$. 
\item[\textbullet] $P_{hwg,t} \in \mathbb{R}_+$: Amount of hot water produced by hot water generator [kW] over the time interval $[t,(t+1)]$. 
\item[\textbullet] $P_{ct,t} \in \mathbb{R}_+$: Amount of condenser water input to the cooling towers [kW] over the time interval $[t,(t+1)]$. 
\item[\textbullet] $P_{cw,t} \in \mathbb{R}$: Net charge/discharge rate [kW] of the chilled water energy storage over the time interval $[t,(t+1)]$. If $P_{cw,t}>0$, the chilled water is being discharged and if $P_{cw,t}<0$ the chilled water is being charged.
\item[\textbullet] $P_{hw,t} \in \mathbb{R}$: Net charge/discharge rate [kW] of the hot water energy storage over the time interval $[t,(t+1)]$. If $P_{hw,t}>0$, the hot water is being discharged and if $P_{hw,t}<0$ the hot water is being charged.
\item[\textbullet] $P_{hx,t} \in \mathbb{R}_+$: Amount of hot water input to the dump heat exchanger (HX) [kW] over the time interval $[t,(t+1)]$. 
\item[\textbullet] $r^{e}_{t}\in \mathbb{R}$: Residual electrical load [kW] over the time interval $[t,(t+1)]$.
\item[\textbullet] $r^{w}_{t}\in \mathbb{R}$: Residual water demand [gal/h] over the time interval $[t,(t+1)]$.
\item[\textbullet] $r^{ng}_{t}\in \mathbb{R}$: Residual natural gas demand [kW] over the time interval $[t,(t+1)]$.
\item[\textbullet] $S^{un}_{cw,t}\in \mathbb{R}$: Slack variable for unmet chilled water load [kW] over the time interval $[t,(t+1)]$.
\item[\textbullet] $S^{ov}_{cw,t}\in \mathbb{R}$: Slack variable for overmet (or over-produced) chilled water load [kW] over the time interval $[t,(t+1)]$.
\item[\textbullet] $S^{un}_{hw,t}\in \mathbb{R}$: Slack variable for unmet hot water load [kW] over the time interval $[t,(t+1)]$.
\item[\textbullet] $S^{ov}_{hw,t}\in \mathbb{R}$: Slack variable for overmet (or over-produced) hot water load [kW] over the time interval $[t,(t+1)]$.
\end{itemize}

\noindent \textbf{States:}
\begin{itemize}
\item[\textbullet] $E_{5,t}\in \mathbb{R}_+$: Energy level of the chilled water energy storage [kWh] at time $t$.
\item[\textbullet] $E_{6,t}\in \mathbb{R}_+$: Energy level of the hot water energy storage [kWh] at time $t$.
\item[\textbullet] $R_{t}$: Peak residual electrical load observed up to time $t$ [kW].
\item[\textbullet] $ul_{cw,t}\in \mathbb{R}$: State variable integrating the slack variable for unmet chilled water load [kW] over the time interval $[t,(t+1)]$.
\item[\textbullet] $ol_{cw,t}\in \mathbb{R}$: State variable integrating the slack variable for overmet (or over-produced) chilled water load [kW] over the time interval $[t,(t+1)]$.
\item[\textbullet] $ul_{hw,t}\in \mathbb{R}$: State variable integrating the slack variable for unmet hot water load [kW] over the time interval $[t,(t+1)]$.
\item[\textbullet] $ol_{hw,t}\in \mathbb{R}$: State variable integrating the slack variable for overmet (or over-produced) hot water load [kW] over the time interval $[t,(t+1)]$.
\end{itemize}

\noindent \textbf{Economic metrics:}
\begin{itemize}
\item[\textbullet] $\Phi^{nocp}(\xi) := {\pi^D}\max_{t\in\mathcal{M}}L_t(\xi)$: Total cost when there is no HVAC central plant in campus [\$].
\item[\textbullet] $\Phi^{sto}(\xi)$: Total cost for stochastic MPC under the realization $\xi$ [\$].
\item[\textbullet] $\Phi^{det}(\xi)$: Total cost for deterministic MPC under the realization $\xi$ [\$].
\item[\textbullet] $\Phi^{perf}(\xi)$: Total cost for perfect information MPC under the realization $\xi$ [\$].
\item[\textbullet] $\textrm{CCP}^{perf}(\xi)$: Ideal cost of central plant under perfect information MPC under realization $\xi$ [\$].
\item[\textbullet] $\textrm{CCP}^{sto}(\xi) $: Cost of central plant under stochastic MPC and realization $\xi$ [\$].
\item[\textbullet] $\textrm{CCP}^{det}(\xi) $: Cost of central plant under deterministic MPC and realization $\xi$ [\$].
\item[\textbullet] $\textrm{VSMPC}(\xi)$: Value of stochastic MPC under the realization $\xi$ [\$].
\end{itemize}

\bibliography{central_plant.bib}

\end{document}